\documentclass{amsart}
\usepackage[marginratio=1:1]{geometry}

\usepackage[hidelinks]{hyperref}
\usepackage{calc}
\newsavebox\CBox
\newcommand\hcancel[2][0.5pt]{%
  \ifmmode\sbox\CBox{$#2$}\else\sbox\CBox{#2}\fi%
  \makebox[0pt][l]{\usebox\CBox}%
  \rule[0.5\ht\CBox-#1/2]{\wd\CBox}{#1}}
\usepackage{blindtext}
\usepackage{color}
\usepackage{hyperref,url}
\usepackage{float}
\usepackage{hyperref}
\usepackage{enumerate}
\usepackage{enumitem}   
\usepackage{bbm}
\usepackage{cancel}
\usepackage{mathrsfs}
\usepackage{colonequals}
\usepackage{pdfpages}

\usepackage{amsmath,amsfonts,amsthm,bm}
\usepackage{mathrsfs}  
\usepackage{xcolor}
\usepackage{amsfonts}
\usepackage{amssymb}
\usepackage{amsmath}
\usepackage{mathtools}
\usepackage{mathrsfs}  

\usepackage[normalem]{ulem}

\numberwithin{equation}{section}
\usepackage[toc,page]{appendix}
\usepackage{tikz-cd}
\theoremstyle{definition}

\newtheorem{definicao}{Definition}[section]
\newtheorem{remark}[definicao]{Remark}

\theoremstyle{plain}

\newtheorem{teorema}[definicao]{Theorem}
\newtheorem{suposicao}{Hypothesis}

\newtheorem{lema}[definicao]{Lemma}
\newtheorem{proposicao}[definicao]{Proposition}
\newtheorem{corolario}[definicao]{Corollary}

\newenvironment{TheoremproofC}[1]{\par\noindent{\textit{Proof of Theorem 3.4.}} \space#1}{\leavevmode\unskip\penalty9999 \hbox{}\nobreak\hfill\quad\hbox{$\qed$}}
\usepackage{scalerel}

\newenvironment{TheoremproofE}[1]{\par\noindent{\emph{Proof of Proposition \ref{conmigoestabachata}.}} \space#1}{\leavevmode\unskip\penalty9999 \hbox{}\nobreak\hfill\quad\hbox{$\qed$}}
\definecolor{roxo}{rgb}{0.44, 0.16, 0.39}
\definecolor{ao(english)}{rgb}{0.0, 0.5, 0.0}
\definecolor{dmagenta}{RGB}{139, 0, 139}
\definecolor{dgreen}{RGB}{0,90,0}
\definecolor{navy}{RGB}{0,0,128}

\usepackage{stackengine}

\def\d{\mathrm d}
\def \d {\mathrm{d}}

\definecolor{iblue}{RGB}{0, 35, 194}
\title[Existence of horseshoes]{Horseshoes for a class  of nonuniformly expanding random dynamical systems on the circle}
\author
{Jeroen S.W. Lamb${}^*$, Giuseppe Tenaglia${}^*$ and Dmitry Turaev${}^*$
}

\allowdisplaybreaks
\begin{document}



\maketitle
\vspace{-1.3cm}
\begin{align*}
   \small {}^* \textit{Imperial College London} 
\end{align*}

\begin{abstract}
{We propose a notion of random horseshoe for one-dimensional random dynamical systems. We prove the abundance of random horseshoes for a class of circle endomorphisms subject to additive noise, large enough to make the Lyapunov exponent positive. In particular, we provide conditions which guarantee that given any pair of disjoint intervals, for almost every noise realization, there exists a positive density sequence of return times to these intervals such that the induced dynamics are the full shift on two symbols.}



\end{abstract}

\section{Introduction and motivation}
One of the most fascinating aspects of dynamical systems is the study of chaos and its properties. In random dynamical systems, the existence of a physical measure \cite{young2002srb} with positive top Lyapunov exponent is a popular notion of chaos. Lyapunov exponents are   statistical quantities  measuring the infinitesimal rate of expansion/contraction along a given direction for typical orbits.
If the top one is positive, then the local expansion eventually separates nearby orbits, thus providing a mechanism  for causing possibly very complicated dynamics. 

In the deterministic setting, a remarkable result by Katok \cite{PMIHES_1980__51__137_0} establishes that  the dynamics of $\mathcal{C}^{1+\alpha}$  diffeomorphisms preserving an ergodic measure with  positive top Lyapunov exponent can be approximated by uniformly hyperbolic invariant sets called horseshoes.\\
For a  diffeomorphism $f$ of a  manifold $X$, two boxes $I_0,I_1 \subset X$, and $N \in \mathbb{N}$, the set 
\begin{align*}
H = \bigcap_{j \in \mathbb{Z}} f^{-jN}(I_0 \cup I_1).
\end{align*}
is called a horseshoe if the differential map $df^N\mid_{T(I_0 \cup I_1)}$ is uniformly hyperbolic  and the map $f^N\mid_{H}$ is topologically conjugated to the full shift on two symbols, in the sense that the bi-infinite sequences of zeros and ones are in one-to-one correspondence with the orbits of points in $H$: 
\begin{equation}\label{simboli}
\forall  \{s_n\}_{n \in \mathbb{Z} } \in \{0,1\}^{\mathbb{Z}},\,\,\exists\,!\, x \in H \mid  f^{Nj}(x) \in I_{s_j} \qquad \forall j \in \mathbb{Z}.
\end{equation}

   The topological conjugacy implies that $f^N\mid_{H}$ has positive topological entropy and is sensitively dependent on initial conditions. The Katok's result connects a  statistical property, the positivity of the top Lyapunov exponent, to topological dynamical  properties. 
   
When aiming to establish a similar correspondence for random dynamics, one encounters a number of difficulties. The main issue is to control return times. In deterministic settings  an essential ingredient in Katok's construction \cite{PMIHES_1980__51__137_0} is the identification of  two boxes ($I_0,I_1$) with periodic returns. In random systems, however, return times typically depend on noise realizations and may be arbitrarily large with positive probability. 

In the present paper, we resolve the problem for a class of circle maps as described in Section $2$. These maps display noise-induced chaos: for small noise amplitude typical trajectories synchronize and converge to a random fixed point, but, beyond a critical noise amplitude,  the dynamics become chaotic. Thus, for moderately large noise, we establish both the positivity of Lyapunov exponent and  abundance of the horseshoe-like dynamics.

Given a set $\Omega$, a map  $\theta \colon \Omega \to \Omega$ and a family of maps $\{f_{\omega}\}_{\omega \in \Omega}$ with $f_{\omega} \in \mathcal{C}^1(\mathbb{S}^1)$,
consider the iterations
\begin{align*}
f^n_{\omega}(x) := f_{\theta^{n-1}(\omega)} \circ \dots \circ f_{\omega}(x).
\end{align*}
\begin{definicao}
Let $I_0,I_1$ be two intervals such that
$I_0 \cap I_1 = \emptyset $. Let $\omega \in \Omega$ and $\kappa >1$. Then the pair $(I_0,I_1)$ is  an $(\omega,\kappa)$-horseshoe if there exists a sequence of  return times $\{n_k(\omega) \}_{k \ge 0}$ satisfying
\begin{equation}\label{positiveDensity}
\limsup_{k \to \infty} \frac{n_k(\omega)}{k} < \infty ,
\end{equation}
and also a family of subsets $\{J(k,\omega)_{i,j} \}_{\{k \ge 0, i,j \in \{0,1\}\}}$  such that  $J(k,\omega)_{i,j} \subset I_i$,
\begin{equation}\label{e1} f^{n_{k+1}(\omega)-n_k(\omega)}_{\theta^{n_k(\omega)}(\omega)}(J(k,\omega)_{i,j}) = I_j,
\end{equation}
and 
\begin{equation}\label{e2}
\left|df^{n_{k+1}(\omega)-n_k(\omega )}_{\theta^{n_k(\omega)}(\omega)}\vert_{J(k,\omega)_{i,j}}\right| >\kappa \qquad \forall  k \ge 0.
\end{equation}
\end{definicao}
In our setting we endow $\Omega$ with a structure of probability space, so that the iterations  $f^n_{\omega}$ represent a random dynamical system. In  Theorem \ref{tea}, which is our main result, we describe a class of random dynamical systems for which we establish the abundance of random horseshoes.

This  result is exemplified by the following special case:

\begin{teorema}\label{Teorema}
Let $\{\omega_n \}_{n \ge 0}$ be i.i.d. random variable uniformly distributed  on $[-\sigma,\sigma] \subset [-\frac{1}{2},\frac{1}{2}]$ and let $L>0$. Consider the random dynamical system (RDS) generated by the iterations $f^n_{\omega}:=f_{\omega_{n-1}} \circ \dots \circ  f_{\omega_0}$, where $f_{\omega_i}$ are circle maps of the form  
\begin{equation}\label{example}
f_{\omega_i}(x) := L\sin(2\pi (x+\omega_i))+ a \,\,\, (\mathrm{mod} \,1),
\end{equation}
with $a \in (0,1)$.
Let $\mathbb{P}:= \left(\frac{\mathrm{Leb}_{[-\sigma, \sigma]}}{2\sigma}\right)^{\mathbb{N}}$   be the natural product measure on $\Omega:= [-\sigma,\sigma]^{\mathbb{N}}$ and let $\theta$ denote the shift on $\Omega$.
Then, there exist constants $\gamma>0$ and $\kappa >1$ such that   for every $a \in (0,1)$, $L\ge 9$ and  $\sigma \ge \frac{11}{10}L^{-\frac{1}{2}}$
\begin{itemize}
    \item[(i)] the Lyapunov exponent is positive:
    \begin{align*}
    \liminf_{n \to \infty}\frac{1}{n}\log |df^n_{\omega}(x)|>\gamma \log(L) \qquad \mathbb{P} \times \mathrm{Leb} \,\,\text{a.s.}
    \end{align*}
    \item[(ii)] there exists a  full $\mathbb{P}$-measure set $\tilde{\Omega} \subset \Omega$  such that for every $\omega \in \tilde{\Omega}$ every pair of connected disjoint intervals $(I_0,I_1)$  is an $(\omega,k)$-horseshoe.
    
\end{itemize}
\end{teorema}

Theorem \ref{Teorema} ensures the existence of horseshoe-like dynamics for a family of non-uniformly expanding RDSs. 
Results similar to Theorem \ref{Teorema} were previously obtained in the setting of uniformly hyperbolic RDSs, either mixing on fibers \cite{huang2016ergodic}, or with an equicontinuous or quasi-periodic base \cite{WenHUANG:281}. Without assuming hyperbolicity, the existence of a weak horseshoe was established for random dynamical systems with positive topological entropy \cite{https://doi.org/10.1002/cpa.21698}. However, weak horseshoes essentially only characterize dynamics with a finite time horizon, whilst our results concern the behaviour for infinite time.

Our general result, Theorem \ref{tea}, 
 holds for random circle endomorphisms  characterized by the following conditions:
\begin{enumerate}
\item\label{A} positive Lyapunov exponent;
    \item\label{B}  sufficiently fast decay of the tail of the annealed distribution of the  Lyapunov exponent; 
    \item\label{C}  uniformly expanding interval that is mapped onto the whole circle;
    \item\label{D} sufficiently large noise amplitude such that there is a positive probability for any given interval to be eventually mapped onto the whole circle.
\end{enumerate}
Then, under further mild conditions we prove that  {\em every pair of disjoint intervals $(I_0,I_1)$ is almost surely an $(\omega,\kappa)$-horseshoe for some $\kappa >1$}.  In particular, these $(\omega,\kappa)$-horseshoes imply chaos in terms of sensitive dependence on initial conditions and positive topological entropy. 

As an application of Theorem \ref{tea}, in Theorem \ref{teb}, we establish horseshoe-like dynamics for predominantly expanding random circle endomorphisms with sufficiently large additive noise.
In particular, in the absence of noise, the RDS in \eqref{example} becomes the deterministic map 
 \begin{align*}
f(x):= L\sin(2\pi x)+ a \,\,(\mathrm{mod}\,1).
\end{align*}
For  typical parameter values, it has sinks and  the complement of their basin of attraction is an uniformly hyperbolic repeller. For this class of examples, Theorem \ref{teb} establishes, for large enough $L$, that the RDS obtained from this endomorphism  by the inclusion of additive noise displays a transition to chaos beyond a specific noise amplitude level, in the sense that the RDS has an ergodic stationary measure with full support, positive Lyapunov exponent, and an abundance of horseshoe-like dynamics. In effect, the noise facilitates the merger of the deterministic attractor and repeller, to the extent that chaotic dynamics prevail. 

This paper is organized as follows:  In Section $2$, we recall the basic setup of random dynamical systems, provide sufficient conditions for the existence of a horseshoe  and state the main result. We also define the class of predominantly expanding random maps. In Section $3$,  we revisit the theory of hyperbolic times \cite{alves2003physical} in application to the random setting. In Section $4$  we obtain  estimates on the quenched tail of hyperbolic times.
In Section $5$ we prove the main Theorem. In section $6$ and $7$ we show that the conclusions of the main theorem hold for the class of predominantly expanding random dynamical systems defined in Section $2$, i.e. Theorem \ref{Teorema} is proved.
\section{Statement of the results}
\subsection{Hypotheses and main result}
Let $(X,\Sigma,\nu)$ be a probability space and denote by $\mathbb{S}^1$ the unit circle, endowed with the arc-length distance  $\mathrm{dist}(\cdot,\cdot)$ and the Borel $\sigma$-algebra $\mathcal{B}(\mathbb{S}^1)$ induced by this distance.  Let also $\mathrm{Leb}(\cdot)$ denote the Lebesgue measure in $\mathbb{S}^1$.

 We consider skew product systems of the form
\begin{equation}\label{skewproduct}
\Theta \colon \Omega \times \mathbb{S}^1 \to   \Omega \times \mathbb{S}^1, \qquad  
\Theta(\omega,x):= (\theta(\omega),g_{\omega}(x)),
\end{equation}
where $\Omega:= X^{\otimes \mathbb{N}}$ is a probability space equipped with the product measure $\mathbb{P} = \nu^{\otimes \mathbb{N}}$ and the product $\sigma$-algebra $\mathcal{F}$ and
\begin{align*}
\theta \colon \Omega \to \Omega, \qquad 
\theta(\{\omega_i\}_{i \ge 0}) = \{\omega_{i+1}\}_{i \ge 0}
\end{align*}
is the shift map acting on $\Omega$. Thus, $(\Omega,\mathcal{F},\mathbb{P},\theta)$ is an ergodic dynamical system. Moreover, $\{g_{\omega} \}_{\omega \in \Omega}$ is a family of circle endomorphisms.

We assume that $\{ g_{\omega}\}_{\omega \in \Omega}$ is a family of continuous  piecewise $\mathcal{C}^2$ endomorphisms such that, for all $\omega \in \Omega$, outside a finite set $\mathcal{S}_{\omega}$ of non-differentiability points, $g_{\omega}$ is $\mathcal{C}^2$ and its first and second derivatives are uniformly bounded. 
Furtermore, define the critical set
\begin{align*}
\mathcal{C}_{\omega} := \{x \in \mathbb{S}^1 \setminus \mathcal{S}_{\omega}\, \,\mid\,\, dg_{\omega}(x)=0 \},
\end{align*}
 and let $\mathcal{SC}_{\omega}:= \mathcal{S}_{\omega} \cup \mathcal{C}_{\omega}$ be the 
 singular set.
 We say that the families $\{S_{\omega}\}_{\omega \in \Omega}$, $\{\mathcal{C}_{\omega}\}_{\omega \in \Omega}$ and $\{\mathcal{SC}_{\omega}\}_{\omega \in \Omega}$ are respectively  the random non-differentiablity set, the random critical set and the random singular set for $\{g_{\omega}\}_{\omega \in \Omega}$. 
The following definition is an adaptation of  \cite[Definition 4.1]{alves2003physical}.
\begin{definicao}\label{cset}
Let $\{g_{\omega}\}_{\omega \in \Omega}$ be as above and let $\{\mathcal{SC}_{\omega}\}_{\omega \in \Omega}$ be its random singular set. We say that $\{\mathcal{SC}_{\omega}\}_{\omega \in \Omega}$ is $\mathbb{P}$ a.s. regular if there exist constants $B>1$ and $\beta > 0$ such that, for almost surely all $\omega \in \Omega$
\begin{itemize}
    \item for every $x \in \mathbb{S}^1\setminus \mathcal{SC}_{\omega}$
    \begin{align*}
      |dg_{\omega}(x)| \ge \frac{1}{B}\mathrm{dist}(x,\mathcal{SC}_{\omega})^{\beta},
    \end{align*}
    \item for every $x,y \in \mathbb{S}^1\setminus \mathcal{SC}_{\omega}$ such that $\mathrm{dist}(x,y) \le \mathrm{dist}(x,\mathcal{SC}_{\omega})/2$
    \begin{align*}
     \left|\log\left|dg_{\omega}(x)\right|- \log\left|dg_{\omega}(y)\right| \right| \le B \frac{\mathrm{dist}(x,y)}{\mathrm{dist}(x,\mathcal{SC}_{\omega})^{\beta}}.
    \end{align*}
\end{itemize}
\end{definicao}
\vspace{1 cm}
The setting of this paper concerns a family $\{g_{\omega}\}_{\omega \in \Omega}$ of circle mappings satisfying the following Hypothesis:
\renewcommand\thesuposicao{(H0)} 
\begin{suposicao}\label{(A)}
The family  $\{g_{\omega}\}_{\omega \in \Omega}$ of  circle endomorphisms satisfies:
 \begin{itemize}
     \item[(i)] $g_{\omega}$ depends only on $\omega_0$, and for all $x \in \mathbb{S}^1$ the maps $\omega \to g_{\omega}(x)$ 
are measurable;
    \item[(ii)] for $\mathbb{P}$ a.s. $\omega \in \Omega$, the maps $x \to g_{\omega}(x)$
    are continuous and piecewise $\mathcal{C}^2$ and  their
     random singular set is regular.  Furthermore
    \begin{align*}
    \operatorname*{ess\,sup}_{\omega \in \Omega} ||g_{\omega}||_{\mathcal{C}^2(\mathbb{S}^1\setminus\mathcal{S}_{\omega})}< \infty.
    \end{align*}
    \end{itemize}
\end{suposicao}

Given a family $\{g_{\omega} \}_{\omega \in \Omega}$ of circle endomorphisms satisfying Hypothesis \ref{(A)}, the random composition $g^n_{\omega}$
\begin{equation}\label{iterazioni}
g^n_{\omega}(x)= g_{\theta^{n-1}(\omega)}(x) \circ \dots \circ g_{\omega}(x),
\end{equation}
 is a Markov process
with associated transition kernel  $\left\{ \mathcal{P}_x(dy):=\mathbb{P}(g_{\omega}(x) \in dy)\right\}_{x \in \mathbb{S}^1}$.
We recall that $\mu$ is a stationary measure for this Markov process if 
\begin{align*}
\mu(A) = \int_{\mathbb{S}^1}\mathcal{P}_x(A) d\mu(x),
\end{align*}
for all Borel sets $A$. Furthermore, $\mu$ is ergodic if and only if the measure $\mathbb{P} \otimes \mu$ is ergodic for the skew product $\Theta$ introduced in \eqref{skewproduct}. If $\mu$ is  ergodic then, by Birkhoff ergodic theorem 
\begin{align*}
\lim_{n \to \infty} \frac{1}{n} \sum_{i=0}^{n-1}f(g^i_{\omega}(x)) = \int_{\mathbb{S}^1} f(x) d\mu(x)     \qquad \mathbb{P}\otimes \mu \,\text{a.s},
\end{align*}
 for all $f \in L^1(\mu)$.   For a more detailed background on Markov processes, see \cite{norris1998markov,hairer2010convergence}.

In order to avoid degenerated behaviour of the random dynamical system, we introduce the following Hypothesis:
\renewcommand\thesuposicao{(H1)} 
\begin{suposicao} \label{(U)} 
The Markov process associated to $g^n_{\omega}$ has an  unique ergodic stationary measure $\mu$ equivalent to the Lebesgue measure. Furthermore, $(\omega,x) \mapsto \log|dg_{\omega}(x)| \in L^1\left(\mu \otimes \mathbb{P}\right)$.
\end{suposicao}
 Note that the ergodic stationary measure $\mu$ is equivalent to Lebesgue if and only if, given any initial condition $x$, the Markov chain $g^n_{\omega}$ starting at $x$ can reach every open set with positive probability (cf. also \cite{blumenthal2022positive,Lian2012PositiveLE}).
 The ergodicity of $\mu$ and the integrability of the function  $\log|dg_{\omega}(x)|$ with respect to $\mathbb{P}\otimes \mu$ is required in \ref{(U)} to ensure the almost sure existence of the Lyapunov exponent \cite{arnold2013random}:
 \begin{definicao}\label{Lyapunov}
Let $g^n_{\omega}$ be the RDS as in \eqref{iterazioni}. The Lyapunov exponent associated to $g^n_{\omega}$ is defined as
\begin{align*}
\lambda_0(\omega,x):= \liminf_{n \to \infty} \frac{1}{n} \log\left|\left(dg^n_{\omega}\right)(x) \right| = \liminf_{n \to\infty}   \frac{1}{n} \sum_{i=0}^{n-1} \log |dg_{\theta^i(\omega)}(g^i_{\omega}(x))|.
\end{align*}
\end{definicao}
\vspace{0.3 cm}
The assumption of ergodicity of $\mu$ in \ref{(U)} implies that the Lyapunov exponent is  $\mathbb{P} \otimes \mu$ almost surely constant and equal to 
\begin{align*}
\lambda_0 = \int_{\Omega}\int_{\mathbb{S}^1}\log|\d g_{\omega}(x)|d\mu(x)d\mathbb{P}(\omega).
\end{align*}

For the class of systems we study, we require that $\lambda_0$ is positive and that we have  estimates  on how much the finite-time ergodic averages deviate from  $\lambda_0$, as stated in the next Hypothesis.
\renewcommand\thesuposicao{(H2)} 
\begin{suposicao} \label{(L)} 
Let $S_n(\omega,x):= \sum_{i=0}^{n-1}\log|dg(g^i_{\omega}(x))|$. Then, there exists $\lambda >0$, $\alpha_1 > 0$ and $C_1 \ge 1$ such that, for  Lebesgue almost all $x \in \mathbb{S}^1$
\begin{align*}
\mathbb{P}\left( S_n(\omega,x)< \lambda n\right) \le \frac{C_1}{n^{4+\alpha_1}}, \qquad \forall n \ge 0.
\end{align*}
\end{suposicao}
Positivity of $\lambda_0$ follows combining Hypotheses  \ref{(U)}-\ref{(L)} with a Borel-Cantelli argument. In particular, one has  
\begin{align*}
\lambda_0 := \lim_{n \to \infty } \frac{|dg^n_{\omega}(x)|}{n} \ge \lambda \qquad \mathbb{P}\otimes \mu\,\,\text{a.s.}
\end{align*}

In order to prove our result, we need further assumptions on the behaviour of the orbits, essentially ensuring that typical orbits do not accumulate too fast around the singular set.  Given $\delta  > 0$ and $\omega \in \Omega$, we define the $\delta$-truncated distance from $x$ to the singular set $\mathcal{SC}_{\omega}$  as 
\begin{equation}\label{criticaldistance}
\mathrm{dist}_{\delta}(x,\mathcal{SC}_{\omega}):= \begin{cases} 1
\qquad &\text{if}\,\, \mathrm{dist}(x,\mathcal{SC}_{\omega})> \delta, \\
\mathrm{dist}(x,\mathcal{SC}_{\omega}) \qquad &\text{otherwise}.
\end{cases}
\end{equation}
\begin{definicao}\label{slow recurrence}
We say that $f^n_{\omega}$ exhibits slow recurrence to the random singular set  $\{\mathcal{SC}_{\omega}\}_{\omega \in \Omega}$, if, for every $\varepsilon>0$, there exists a $\delta> 0$ such that, for $ \mathbb{P} \otimes \mu $ almost all $(\omega,x)  \in \Omega \times \mathbb{S}^1$
\begin{align*}
\limsup_{n \to \infty} \frac{1}{n} \sum_{j=0}^{n-1}-\log\mathrm{dist}_{\delta}(g^j_{\omega}(x),\mathcal{SC}_{\theta^{j}(\omega)})< \varepsilon.
\end{align*}
\end{definicao}
\vspace{0.8 cm}
As a consequence of the ergodicity assumption \ref{(U)} we have 
\begin{align*}
\lim_{n \to \infty} \frac{1}{n} \sum_{j=0}^{n-1}-\log\mathrm{dist}_{\delta}(g^j_{\omega}(x),\mathcal{SC}_{\theta^{j}(\omega)}) = \int_{\Omega} \int_{\mathbb{S}^1}-\log\mathrm{dist}_{\delta}(x,\mathcal{SC}_{\omega}) d\mu(x)d\mathbb{P}(\omega) \qquad  \mathbb{P}  \otimes \mu  \,\, \text{a.s.}
\end{align*}
 In order to control the deviation in the convergence of the above limit, so to ensure that typical orbits do not get close to the singular set too often,  we introduce 
\renewcommand\thesuposicao{(H3)} 
\begin{suposicao} \label{(M)}
Let $Z_n(\omega,x,\delta):= \sum_{j=0}^{n-1}\log(\mathrm{dist}_{\delta}(g^j_{\omega}(x),\mathcal{SC}_{\theta^j(\omega)}))$. There exists a function $H(\delta)$ that converges to zero as $\delta \to 0$ and, constants $\alpha_2 >1$ and $C_2 \ge 1$ such that, for Lebesgue almost all $x \in \mathbb{S}^1$ and all $\delta>0$ suffciently small
\begin{align*}
\mathbb{P}\left( Z_n(\omega,x,\delta)> H(\delta)n \right) \le \frac{C_2}{n^{4+\alpha_2}}.
\end{align*}
\end{suposicao}
Combining Hypotheses  \ref{(U)} and \ref{(M)}, we see that  there exists $\tilde{H}=\tilde{H
}(\delta)\le H(\delta)$ such that $\tilde{H} \to 0$ as $\delta \to 0$ and
\begin{equation}\label{slowslow}
\lim_{n \to \infty} \frac{1}{n} \sum_{j=0}^{n-1}-\log\mathrm{dist}_{\delta}(g^j_{\omega}(x),\mathcal{SC}_{\theta^j(\omega)}) \le  \tilde{H}(\delta)  \qquad  \mathbb{P}\otimes \mu\,\,\text{a.s.}
\end{equation}
In Section $7$, we show that Hypotheses \ref{(L)}-\ref{(M)} are satisfied by a class of non-uniformly expanding RDSs with sufficiently large noise amplitude. \\
In order to state the next Hypothesis, we define the concept of full branch expanding intervals, which will be crucial to the horseshoe construction.
\begin{definicao}\label{fullBranchExpanding}
Let $I \subset \mathbb{S}^1$ be an interval. For $\omega \in \Omega$ and $\kappa_0>1$, we say that $I$ is $(\omega,\kappa_0,K)$-\emph{full branch expanding}  if there exists $J \subset I$ and $i \le K$ such that
$g^i_{\omega}(J) = \mathbb{S}^1$ and $|dg^i_{\omega}\mid_{J}|> \kappa_0$.
\end{definicao}
In Hypothesis \ref{(O)}, we require that any  interval  is  $(\omega,\kappa_0,K)$-full branch expanding with positive probability. This probabilistic assumption on the finite-time growth  of intervals, combined with the asymptotic behaviour given in Hypotheses \ref{(L)}-\ref{(M)},  ensures that intervals of any size are eventually onto with full probability. This grants us the necessary information to construct the desired horseshoes.
\renewcommand\thesuposicao{(H4)} 
\begin{suposicao} \label{(O)}
There exists $\kappa_0>1$ such that, for all sufficiently small $\eta>0$, there exist a natural number $K$ and  a real number $\iota \in (0,1)$ such that, for any interval $I \subset \mathbb{S}^1$ satisfying $\frac{\eta}{4}\le |I| \le 4\eta$, we have
\begin{equation}\label{ctrol}
\mathbb{P}\left\{I\, \text{is  $(\omega,\kappa_0,K)$-full branch expanding}\right\} > \iota.
\end{equation}
\end{suposicao}
In Hypothesis \ref{(P)} we additionally require that for all $\omega \in \Omega$ there exists an interval $I_{\omega}$ which is full branch and uniformly expanding for $g_{\omega}$. This assumption is relevant for two reasons. First, If $\{g_{\omega}\}_{\omega \in \Omega}$ satisfies \ref{(P)}, then any $(\omega,\kappa_0,n)$-full branch expanding  is also  $(\omega,\kappa_0,m)$-full expanding  branch  for all $m \ge n$. 
This is a key property used in the proof of Theorem \ref{tea}. Moreover, in Section $7$ we show that \ref{(P)} implies \ref{(O)} in the case of additive noise with sufficiently large noise amplitude.
\renewcommand\thesuposicao{(H5)}
\begin{suposicao} \label{(P)}
For $\mathbb{P}$ almost all $\omega \in \Omega$ there exists an interval $I_{\omega}$ such that $g_{\omega}|_{I_{\omega}}= \mathbb{S}^1$ and 
\begin{equation}
\inf_{\omega \in \Omega}|dg_{\omega}|_{I_{\omega}}|>1.
\end{equation}
\end{suposicao}
In order to state the main result of this paper 
it remains to define  what it means for two intervals $I_0,$ and $I_1$ to form a horseshoe.
\begin{definicao}\label{Horseshoe}
Let $I_0,I_1$ be two intervals such that
$I_0 \cap I_1 = \emptyset $. Let $\omega \in \Omega$ and $\kappa >1$. Then the pair $(I_0,I_1)$ is  an $(\omega,\kappa)$-horseshoe if there exists a sequence of  return times $\{n_k(\omega) \}_{k \ge 0}$ satisfying
\begin{equation}\label{positiveDensity}
\limsup_{k \to \infty} \frac{n_k(\omega)}{k} < \infty ,
\end{equation}
and also a family of subsets  $\{J(k,\omega)_{i,j} \}_{\{k \ge 0, i,j \in \{0,1\}\}}$  such that  $J(k,\omega)_{i,j} \subset I_i$,
\begin{equation}\label{e1} g^{n_{k+1}(\omega)-n_k(\omega)}_{\theta^{n_k(\omega)}(\omega)}(J(k,\omega)_{i,j}) = I_j,
\end{equation}
and 
\begin{equation}\label{e2}
\left|dg^{n_{k+1}(\omega)-n_k(\omega )}_{\theta^{n_k(\omega)}(\omega)}\vert_{J(k,\omega)_{i,j}}\right| >\kappa \qquad \forall  k \ge 0.
\end{equation}
\end{definicao}
\begin{teorema}\label{tea}
Let $g^n_{\omega}$ be an RDS satisfying \ref{(A)}-\ref{(P)}. Then, there exists a full $\mathbb{P}$ measure set $\tilde{\Omega}$ and $\kappa>1$ such that every disjoint pair of intervals $(I_0,I_1)$  is an $(\omega,\kappa)$-horseshoe for all $\omega \in \tilde{\Omega}$. Furthermore, the random times $\{n_k(\omega)\}_{k \ge 0}$ are measurable and  satisfy 
\begin{align*}
\lim_{k \to \infty} \frac{n_k(\omega)}{k} = \mathbb{E}[n_0]\qquad \mathbb{P}\,\text{a.s}.
\end{align*}
\end{teorema}
An immediate corollary is the following:
\begin{corolario}
Let $g^n_{\omega}$ be an RDS satisfying \ref{(A)}-\ref{(P)}. Then, the system displays sensitive dependence on initial conditions. Furthermore, each of the $(\omega,\kappa)$-horseshoes carries positive topological entropy greater or equal to $\frac{\log(2)}{\mathbb{E}[n_0]}$.
\end{corolario}
\subsection{A sketch of the proof of Theorem \ref{tea}}
We proceed to sketch the strategy of the proof of Theorem \ref{tea}, the full details of which are laid out in Section $5$.

We establish the following fact (see Lemma \ref{return times}):
\begin{itemize}
    \item There exists $\kappa>1$ and $\alpha >0$ such that,  for every interval $I$,  the stopping time  
    \begin{align*}
    m(\omega,I):= \min\{ n \ge 0 \mid \exists  J \subset I\,\,\text{such that}\,\, f^n_{\omega}(J) = \mathbb{S}^1\,\, \text{and} \,\,|df^n_{\omega}\,\vert_{J}|> \kappa \} 
    \end{align*}
    satisfies 
    \begin{equation}\label{boundk}
    \mathbb{P}\left(m(\omega,I)> l \right)\le  \frac{D}{l^{3+\alpha}} \qquad   \forall\, l \ge 1,
    \end{equation}
    for some $D \ge 1$ depending on $I$.
\end{itemize}

The proof of this fact proceeds as follows:
\begin{itemize}
    \item First we prove (Proposition \ref{hyperbolic ball} and Theorem \ref{hyperbolictimes})  that there exists $\kappa_1>1$ and a $\delta_1>0$ small enough such that for every interval $I \subset \mathbb{S}^1$  there exists a sequence of almost surely finite stopping times  $\tau_k$ and intervals $J_k$ such that  
    \begin{align*}
    |f^{\tau_k}_{\omega}(J_k) | &= 2\delta_1 \\
    |df^{\tau_k}_{\omega}\,\vert_{J_k} | &> \kappa_1.
    \end{align*}
    This is a classical result in the theory of hyperbolic times \cite{alves2003hyperbolic}.
    \item Since $|f^{\tau_k}_{\omega}(J_k)|=2\delta_1$, Hypothesis \ref{(O)} implies the existence of $K \in \mathbb{N}$, $\kappa_0 >1$ and  and $\iota \in (0,1)$, both depending on $\delta_1$, such that 
    \begin{equation}\label{a}
    \mathbb{P}\left(  \mathcal{S}'\left(K,f^{\tau_k}_{\omega}(J_k),\kappa_0\right)\right) >i,
    \end{equation}
    where we denoted
    \begin{align*}
     \mathcal{S}'\left(K,f^{\tau_k}_{\omega}(J_k),\kappa_0 \right):= \left\{\omega \in \Omega \colon f^{\tau_k}_{\omega}(J_k)\, \text{is  $(\theta^{\tau_k}(\omega),\kappa_0,K)$-full branch expanding}\right\}.
    \end{align*}
    Using \eqref{a} and the independence property of the product measure $\mathbb{P}$ we prove Proposition \ref{almostDone}, which asserts that
    \begin{align*}
    \mathbb{P}\left( \bigcap_{k \ge 0} \mathcal{S}'\left(K,f^{\tau_k}_{\omega}(J_k),\kappa_0\right)^c\right) =0,
    \end{align*}
   where  $\mathcal{S}'\left(K,f^{\tau_k}_{\omega}(J_k),\kappa_0\right)^c$ is the complement of the event $\mathcal{S}'\left(K,f^{\tau_k}_{\omega}(J_k),\kappa_0 \right)$ in $\Omega$.
    As a result, there exists $\mathbb{P}$ almost surely a $m_0(\omega)$, an integer $i \le K$ and an interval $J'(\omega) \subset J_{m_0(\omega)} \subset I$ such that 
    \begin{align*}
    f^{\tau_{m_0}+i}_{\omega}(J')&=\mathbb{S}^1  
    \end{align*}
    and 
    \begin{align*}
    |df^{\tau_{m_0}+i}_{\omega}\mid_{J'}| > \kappa:= \min\{\kappa_0, \kappa_1\}.
    \end{align*}
    This proves the claim, as it establishes the almost sure bound  $m(\omega,I)\le \tau_{m_0}+K$ for the stopping time. The estimate in \eqref{boundk} is proved in Lemma \ref{return times}.
\end{itemize}
The bounds in \eqref{boundk} are  used in Proposition \ref{almostDone}, assuming Hypothesis \ref{(P)}, to prove the almost sure existence for any two given intervals $I_0,I_1$ of a sequence of horseshoe return times $\{n_k(\omega) \}_{k \ge 0}$. The definition of the random sequence $\{n_k(\omega)\}_{k \ge 0}$ implies that the sequences of increments $\{n_{k+1}(\omega)-n_k(\omega) \}_{k \ge 1}$ and $n_0(\omega)$ are independent and identically distributed. Therefore, Proposition \ref{PosD} establishes positive density (see equation \eqref{positiveDensity} in Definition \ref{Horseshoe}), using the law of large numbers.
Repeating the above construction for all possible pairs of dyadic intervals, we obtain a full $\mathbb{P}$-measure set $\tilde{\Omega}$  such that any pair of disjoint intervals is an $(\omega,\kappa)$-horseshoe.
This concludes the proof of Theorem \ref{tea}.
\subsection{Application to random systems with predominant expansion and large additive noise}


We proceed to construct the class of RDS to which we apply Theorem \ref{tea}. Consider the probability space $(\Omega,\mathcal{F},\mathbb{P})$, where  $\Omega:= [-\sigma,\sigma]^{\mathbb{N}}$, $\mathbb{P}:= \left(\frac{\mathrm{Leb}\vert_{[-\sigma,\sigma]}}{2\sigma }\right)^{\otimes \mathbb{N}}$ and $\mathcal{F}$ is the  product Borel $\sigma$-algebra obtained via Kolmogorov extension theorem. Let $f \colon \mathbb{S}^1 \to \mathbb{S}^1$ be a continuous piecewise $\mathcal{C}^2$ map. Then, similarly to what we have done before, we associate to $f$ the set  $\mathcal{S}$ of non differentiability points of $f$, the critical set $\mathcal{C}:= \{x \in \mathcal{S}^1 \setminus \mathcal{S} \colon df = 0 \}$ and the singular set $\mathcal{SC}:= \mathcal{S}\cup \mathcal{C}$.
\begin{definicao}\label{admissibleRDS}
We say that a RDS $f^n_{\omega}$ is $(\sigma,R)$-predominantly expanding if the maps belonging to the family $\{ f_{\omega}\}_{\omega \in \Omega}$ are of the form $f_{\omega} \colon \mathbb{S}^1 \to \mathbb{S}^1 $
\begin{align*}
 f_{\omega}(x):=f(x+\omega_0),
\end{align*}
where
\begin{enumerate}
\item\label{one} the generating map $f$ satisfies the following conditions:
\begin{itemize}
\item[(i)]\label{uno} $f$ is a piecewise $\mathcal{C}^2$  smooth endomorphism  with bounded first and second derivatives;
\item[(ii)]\label{due}
Near each point of the singular set $\mathcal{SC}$ of $f$, either the first or the second derivative of $f$ is bounded away from zero. 
\item [(iii)]\label{tre}Let $ G := \{x \in \mathbb{S}^1 \setminus \mathcal{S} \mid |df| > R \}$ be the super-expanding region for $f$. Then $G$ is the union of $k \ge 1$ connected components $G_1,\dots,G_k$ of $G$ whose size is  larger than $\frac{1}{R}$. 
\end{itemize}
\item \label{two}
The following inequality holds:
\begin{equation}\label{reF}
\frac{1}{\log(R)}\int_{C}\log| df(x)|dx  + \frac{2}{R}+D(R)
>0.
\end{equation}
where $C$ is the contracting region of $f$.
\item\label{three}  
The noise is sufficiently large, i.e.
\begin{equation}\label{noiseRange}
\sigma> \frac{1}{R}+D(R),
\end{equation}
where
\begin{equation}\label{dierre}
D(R) :=  \mathrm{Leb}\{x \in \mathbb{S}^1 \colon |df| <R \}= \mathrm{Leb}\left(\mathbb{S}^1 \setminus G\right)
\end{equation}
\end{enumerate}
\end{definicao}

Note that the $(\sigma,R)$-predominant expansion is a generalization of the concept of predominant expansion introduced in \cite{blumenthal2022positive,Lian2012PositiveLE}. 

Condition \eqref{noiseRange} implies that the noise amplitude is large enough to escape the any connected component of $\mathbb{S}^1 \setminus G$ in one iterate, whilst  \eqref{reF} can be interpreted as an estimate for the Lyapunov exponent.

A consequence of \eqref{noiseRange} and \eqref{reF}, that is crucial in the proof of Theorem \ref{teb}, is that there exists $h>0$ satisfying 
\begin{equation}\label{himpli1}
\sigma > \frac{D(R)}{2h}
\end{equation}
and 
\begin{equation}\label{himpl2}
Z(h):= \log(R)(1-h)+ \frac{1}{2\sigma}\int_C \log|df(x)|dx>0.
\end{equation}

 As shown in section $7$, \eqref{himpli1} implies that the asymptotic frequency of visits to $\mathbb{S}^1 \setminus G$ is strictly smaller than $h$. Therefore, we can interpret the parameter $h$ as an upper  bound for the asymptotic frequency of visits to the complement of the super-expanding region $G$.
 
In order to avoid unnecessary complications due do the fact that for $\sigma > \frac{1}{2}$ the family of maps $\omega_0 \to f(x)+\omega_0$ are not injected anymore, we restrict our attention to the case $\sigma< \frac{1}{2}$. We do expect the result to hold true also for $\sigma>\frac{1}{2}$ for our class of maps. Note that by \eqref{noiseRange} this implies also $R>2$.
 
In such settings, for the class of $(\sigma,R)$-predominantly expanding random systems, we apply Theorem \ref{tea} and prove the following result.  
\begin{teorema}\label{teb}
Any $(\sigma,R)$-predominantly expanding RDS satisfies all the Hypotheses \ref{(A)}-\ref{(P)} of Theorem \ref{tea}. As a consequence, there exists $\kappa>1$ such that any pair of disjoint intervals is an $(\omega,\kappa)$-horseshoe for $\mathbb{P}$-almost every $\omega \in \Omega$.
\end{teorema}
As the RDS generated by the family of maps in \eqref{example}, for any $L \ge 3$, is $(\sigma,R)$-predominantly expanding for $R=\min\{3,L^{\frac{1}{2}}\}$ and   $\sigma$ satisfying \eqref{noiseRange}-for a proof, see section $8$- then Theorem \ref{Teorema} is a corollary of Theorem \ref{teb}.
\section{Random hyperbolic times}
\subsection{Definition and key property}
The aim of this section is to develop the essential tool in our analysis: the hyperbolic times. A theory of hyperbolic times for non-uniformly expanding maps was developed by Alves \cite{Alves2000SRBMF,ALVES20001} to construct SRB measures for Viana maps, and then hyperbolic times have been used extensively later to construct Young towers \cite{ASENS_2002_4_35_1_77_0}. In  this section, we extend the notion of hyperbolic times to the setting of RDS, building on techniques from \cite{alves2003physical,ASENS_2002_4_35_1_77_0}. In order to keep the exposition self-contained, we provide complete proofs.  

Let $g^n_{\omega}$ be the RDS induced by \eqref{iterazioni} satisfying \ref{(A)}-\ref{(M)}, therefore there exists $\beta>0$ such that the random critical set of $g^n_{\omega}$ satisfies  Definition \ref{cset}. We fix once and for all some positive number $b<\min\left\{ \frac{1}{2},\frac{\beta^{-1}}{2}\right\}$.  The following definition is adapted from  \cite[Definition 4.5]{alves2003physical}.


    
\begin{definicao}\label{hyperbolic times}
 Given $\kappa_1>1$ and $\delta>0$, we say that $n$ is a $(\kappa_1,\delta)$-hyperbolic time for $(\omega,x) \in \Omega \times \mathbb{S}^1$ if
\begin{align*}
\left|\left(dg^{m}_{\theta^{n-m}(\omega)}\right)(g^{n-m}_{\omega}(x))\right| \ge \kappa_1^{n-m},\qquad \mathrm{dist}_{\delta}(g^{m}_{\omega}(x),\mathcal{SC}_{\theta^{n-m}(\omega)}) \ge \kappa_1^{-bm}
\end{align*}
for all $0 < m \le  n$. 
\end{definicao}

For $x \in \mathbb{S}^1$, let $B(x,r)$ be the open ball around $x$ of radius $r$. 
The following Proposition is an adaptation of \cite[Lemma 4.6]{alves2003physical} and shows that, whenever $n$ is a hyperbolic time for $(\omega,x)$, $g^n_{\omega}$  uniformly expands a sufficiently small neighbourhood of $x$ to a neighbourhood of $g^n_{\omega}(x)$ of a fixed  radius $\delta_1$. This is the essential property of hyperbolic times needed to prove  Theorem \ref{tea}. 
\begin{proposicao}\label{hyperbolic ball}
Suppose $n$ is a $(\kappa,\delta)$-hyperbolic time for $(\omega,x)$. Then, for any $\delta_1>0$ small enough, there exists  an interval $J \subset B\left(x,\delta_1\kappa_1^{-\frac{n}{2}}\right)$ containing $x$ such that $g^n_{\omega}(J)= B(g^n_{\omega}(x),\delta_1)$, $|dg^n_{\omega}\vert_{J}|>\kappa_1^{\frac{n}{2}}$, and  $g^i_{\omega}(J) \subset B\left(x,\delta_1\kappa_1^{-\frac{n+i}{2}}\right)$  for all $i= 1,\dots,n-1$ .
\end{proposicao}
\begin{proof}
First, we state, without proof, a distortion estimate, which is an easy consequence of  \cite[Lemma 4.6]{alves2003physical} and the regularity of the random singular set.
\begin{lema}\label{rdist}
Suppose the random singular set is regular. Given $\kappa_1>1$ and $\delta>0$, if $n$ is a $(\kappa_1 ,\delta)$ hyperbolic time for $(\omega,x)$, then
\begin{align*}
|dg_{\omega}(y)| \ge {\kappa_1}^{-\frac{1}{2}}|dg_{\omega}(x)|,\qquad \forall y \in B\left(x,\delta_1\kappa_1^{-\frac{n}{2}}\right),
\end{align*}
for a sufficiently small $\delta_1>0$.
\end{lema}
We proceed to prove Proposition \ref{hyperbolic ball} by induction on $n$. Let $n=1$ be a $(\kappa_1,\delta)$-random hyperbolic time for $(\omega,x)$.  By Lemma \ref{rdist}, the ball $B\left(x,\delta_1\kappa_1^{-\frac{1}{2}}\right)$ is expanded by a factor $\left(\kappa_1\right)^{\frac{1}{2}}$. Therefore, Proposition \ref{hyperbolic ball} holds for $n=1$. Assuming that the claim holds for $n$, we prove that the claim holds also for $n+1$. 

Let $n+1$ be a $(\kappa_1,\delta)$ hyperbolic time for $(\omega,x)$.
Then, $n$ is a  $(\kappa_1,\delta)$-random hyperbolic time for $(\theta(\omega),g_{\omega}(x))$.
Using the inductive step, there is an interval $I_n \subset B\left(g_{\omega}(x),\delta_1\kappa_1^{\frac{-n}{2}}\right)$ containing $x$  such that
\begin{equation}\label{bypass1}
g^n_{\theta(\omega)}(I_n) = B(g^{n+1}_{\omega}(x),\delta_1)
\end{equation}
and
\begin{equation}\label{bypass2}
g^i_{\theta(\omega)}(I_n) \subset B\left(x,\delta_1\kappa_1^{\frac{-n+i}{2}}\right) \,\,\forall i = 1,\dots, n-1. 
\end{equation}
We claim that there is no other interval  $I \subset B\left(g_{\omega}(x),\delta_1\kappa_1^{-\frac{n}{2}}\right)$, satisfying \eqref{bypass1} and \eqref{bypass2}. Indeed, by definition of  hyperbolic times, we see that either $\mathrm{dist}(g^k_{\omega}(x),\mathcal{SC}_{\theta_k(\omega)}) \ge \delta > \delta_1$ or $\mathrm{dist}(g^k_{\omega}(x),\mathcal{SC}_{\theta^k(\omega)})> \kappa_1^{-b(n+1-k)}> \kappa_1^{-\frac{(n+1-k)}{2}}$. In both cases, $\mathcal{SC}_{\theta^k(\omega)}$ lies outside of $B\left(g^k_{\omega}(x),\delta_1\kappa_1^{\frac{-n-1+k}{2}}\right)$, for all $k = 1,\dots,n$. 
This implies that $g_{\theta^k(\omega)}$ maps ${B\left(g^k_{\omega}(x),\delta_1\kappa_1^{\frac{-n-1+k}{2}}\right)}$ diffeomorphically onto its image.  
In particular, by \eqref{bypass2}, $g^{n}_{\theta(\omega)}|_{I_n}$ is injective. Hence, by \eqref{bypass1}, $I=I_n$ and the claim is proved.

Next, we claim that $g_{\omega}\left(B\left(x,\delta_1\kappa_1^{\frac{-n-1}{2}}\right)\right)$ contains $I_n$. We argue by contradiction. Suppose that  $g_{\omega}\left(B\left(x,\delta_1\kappa_1^{\frac{-n-1}{2}}\right)\right)$ does not contain $I_n$.  Then, at least one of the intervals  $g_{\omega}\left(\left[x-\delta_1\kappa_1^{\frac{-n-1}{2}},x\right]\right)$ and $g_{\omega}\left(\left[x,x+\delta_1\kappa_1^{\frac{-n-1}{2}}\right]\right)$  is strictly contained in $I_n$. 
Assume that $g_{\omega}\left(\left[x-\delta_1\kappa_1^{\frac{-n-1}{2}},x\right]\right) \Subset  I_n$. The other case can be treated analogously. Then, by \eqref{bypass2},
\begin{align*}
g^i_{\omega}\left(\left[x-\delta_1\kappa_1^{\frac{-n-1}{2}},x\right]\right) \Subset B\left(x,\delta_1\kappa_1^{-\frac{n+i}{2}}\right) \qquad \forall i =1,\dots n.
\end{align*}
and, by \eqref{bypass1}
\begin{equation}\label{keyeq}
g^{n+1}_{\omega}\left(\left[x-\delta_1\kappa_1^{\frac{-n-1}{2}},x\right]\right) \Subset B(g^{n+1}_{\omega}(x),\delta_1)
\end{equation}
Using equation \eqref{bypass2} and Lemma \ref{rdist}, we have
\begin{equation}\label{tobeused}
\left|\left(dg^{n+1}_{\omega}\right)(y)\right|= \prod_{i} \left|dg_{\theta^i(\omega)}(g^i_{\omega}(y))\right|
\ge \kappa_1^{-\frac{n+1}{2}}|g^{n+1}_{\omega}(x)|\ge \kappa_1^{\frac{n+1}{2}} \qquad \forall y \in \left[x-\delta_1\kappa_1^{\frac{-n-1}{2}},x\right].
\end{equation} 
As a consequence, 
\begin{align*}
\left|g^{n+1}_{\omega}\left(\left[x-\delta_1\kappa_1^{\frac{-n-1}{2}},x\right]\right) \right|\ge \delta_1,
\end{align*}
which contradicts \eqref{keyeq}.

One concludes the proof of Proposition \ref{hyperbolic ball} by choosing the interval
\begin{align*}
J := g_{\omega}^{-1}(I_n) \cap \left(B\left(x,\delta_1\kappa_1^{\frac{-n-1}{2}}\right)\right).
\end{align*}
Indeed, $x \in J$ by the injectivity of $g_{\omega}$, the map $g^{n+1}_{\omega}$ takes $J$
onto $B(g^{n+1}_{\omega}(x),\delta_1)$  by \eqref{bypass1}, and $\left|dg^{n+1}_{\omega}\vert_{J}\right|> \kappa_1^{\frac{n+1}{2}}$ from computations similar to  \eqref{tobeused}.
\end{proof}
\subsection{Positive frequency of random hyperbolic times}
 The following proposition implies the existence and positive asymptotic frequency for random hyperbolic times. Indeed, we show  that for any $(\omega,x)$, if the orbit at time $N$ exhibits a positive finite-time Lyapunov exponent and a sufficiently slow recurrence to the singular set. Equivalently, using the notation in \ref{(L)}-\ref{(M)}, if $S_N(\omega,x)>\lambda N$ and $Z_N(\omega,x,\delta)<H(\delta)N$,   then the proportion of the hyperbolic times in $\{1,\dots,N\}$ is bounded away from $0$, uniformly for all $(\omega,x)$ and $n$.
This result is used in Section $4$ to obtain estimates on the distribution of the hyperbolic times. 

The statement and the proof of the following proposition are an adaptation of  \cite[Proposition 4.9 ]{alves2003physical}.

\begin{proposicao}\label{hyperbolictimes}

 Given $N \in \mathbb{N}$, let $S_N(\omega,x)=\sum_{i=0}^{N-1}\log\left(|dg_{\theta^i(\omega)}(g^i_{\omega}(x))|\right)$ and ${Z_N(\omega,x,\delta)=\sum_{i=0}^{N-1}-\log\left(\mathrm{dist}_{\delta}(g^i_{\omega},\mathcal{SC}_{\theta^i(\omega)})\right)}$. Let $\lambda $ be as in Hypothesis \ref{(L)} and $H(\delta)$ as in \ref{(M)}. Then, there exist $\kappa_1>1$, $\delta>0$ and $\gamma \in (0,1)$, which are independent of $(\omega,x)$, such that, if $S_N(\omega,x)> \lambda N$ and $Z_N(\omega,x,\delta)< H(\delta)N$,  the number of  $\left(\kappa_1,\delta\right)$-hyperbolic times  for $(\omega,x)$ that lie in $\{1,\dots,N\}$   is larger than $\gamma N$.
\end{proposicao}
The following lemma, due to Pliss,  is a key ingredient is the proof of this proposition. For a reference, see \cite{pliss1972conjecture} or \cite[Lemma 3.9]{alves2003physical}.
\begin{lema}\label{Pliss}
Given $0 < c \le A$ let $\gamma:= \frac{c}{A}$. Assume that $a_1,\dots,a_N$ satisfy $a_j \le A$ for all $1 \le j \le N$ and
\begin{align*}
\sum_{j=1}^{N} a_j \ge cN
\end{align*}
Then, there exist $l \ge \gamma N$ and $1 \le n_1 < \dots < n_l \le N $ such that 
\begin{align*}
\sum_{j=k}^{n_i} a_j  \ge 0
\end{align*}
for every  $k= 1,\dots,n_i$ and $ i=1,\dots,l$.
\end{lema}
\begin{TheoremproofC}
Consider the sequence
\begin{align*}
a_{j+1}:= \log\left|dg_{\theta^j(\omega)}(g^j_{\omega}(x))\right| - \frac{\lambda}{2}, \qquad j=0,\dots N-1
\end{align*}
Due to (ii) in Hypothesis \ref{(A)}, $a_j \le A:= \sup_{\omega}||dg_{\omega}||_{\mathcal{C}^1(\mathbb{S}^1\setminus \mathcal{S}_{\omega})}$ for all $j= 1,\dots,N$.
Since $S_N(\omega,x)>\lambda N$, we have
\begin{align*}
\sum_{j=1}^{N} a_j \ge \frac{\lambda}{2}N.
\end{align*}
By Lemma \ref{Pliss}, then, there exists a sequence of $l_1 \ge \gamma_1 N$ numbers $1\le n_1,\dots,n_{l_1}\le N$ such that
\begin{align*}
\sum_{j=k}^{n_i}a_j \ge 0,
\end{align*}
for all $k = 1,\dots,n_i$ and $i = 1,\dots,l_1$, with $\gamma_1:= \frac{\lambda}{2A}$. Using the definition of $a_j$, we may rewrite the above inequality as
\begin{equation}\label{rectimes}
\left|dg^{n_i-m}_{\theta^m(\omega)}(g^m_{\omega}(x))\right| \ge e^{\frac{\lambda}{2}(n_i-m)},
\end{equation}
for every $m =0,\dots,n_i-1$ and $i = 1,\dots,l$.

Next, we claim that, if  $\delta>0$ is  small enough,  then there exists $\gamma_2 \in (1-\gamma_1,1)$ and  $l_2 \ge \gamma_2 N $ numbers $1 \le \tilde{n}_{1},\dots,\tilde{n}_{l_2}\le N$ satisfying
\begin{equation}\label{aiutami}
\sum_{j=k}^{\tilde{n}_{i}}\log(\mathrm{dist}_{\delta}(g^{j}_{\omega}(x),\mathcal{SC}_{\theta^j(\omega)})) \ge - \frac{b\lambda}{8}(\tilde{n}_{i}-k),
\end{equation}
for all $k = 1,\dots,\tilde{n}_i-1$ and $ i = 1,\dots,l_1$, with $b$ as in Definition \ref{hyperbolic times}.

Indeed, since $H(\delta) \to 0$ as $\delta \to 0$, it is possible to choose $\delta>0$  small enough such that $H(\delta) < \frac{\gamma_1 b \lambda}{8}$. Consider the sequence
\begin{equation}\label{numbers}
\tilde{a}_{j+1}:= \log\mathrm{dist}_{\delta}\left(g^{j}_{\omega}(x),\mathcal{SC}_{\theta^j(\omega)}\right)+ \frac{b\lambda}{8},\qquad j=0,\dots, N-1.
\end{equation}
As $Z_N(\omega,x,\delta)\le H(\delta)N$, we have
\begin{align*}
\sum_{j=1}^N \tilde{a}_j \ge \left( \frac{b\lambda}{8}- H(\delta)\right)N.
\end{align*}
Since $\mathrm{dist}_{\delta}(\cdot,\cdot)\le 1$ by definition, it follows that $\tilde{a}_j \le  \frac{b\lambda}{8}$, for all $j=1\dots,N$.
By Lemma \ref{Pliss}, there exists $l_2 \ge \gamma_2 N $ numbers $1 \le \tilde{n}_1,\dots,\tilde{n}_{l_2}\le N$ such that
\begin{align*}
\sum_{j=k}^{\tilde{n}_i} \tilde{a}_j \ge 0,
\end{align*}
for all $k = 1,\dots,\tilde{n}_i$ and $i=1,\dots,l_2$, where $\gamma_2= 1- \frac{8H(\delta)}{\lambda}$. From the  definition of $\tilde{a}_j$, in \eqref{numbers}   we obtain \eqref{aiutami} as claimed.

To finish the proof, note that $\gamma:= \gamma_1+\gamma_2-1>0$. Indeed
\begin{align*}
\gamma_1+\gamma_2= \gamma_1 + 1-  \frac{8H(\delta)}{\lambda} > \gamma_1 + 1 - \frac{8b\lambda\gamma_1}{8b\lambda} > 1.
\end{align*}
 Since $\left|\{n_1,\dots,n_{l_1}\} \right|\ge \gamma_1N$ and $\left|\{\tilde{n}_1,\dots,\tilde{n}_{l_2}\} \right| \ge \gamma_2N$ and $\gamma:= \gamma_1+\gamma_2-1 >0$, the pigeonhole principle implies that the set  $\{n_1,\dots,n_{l_1}\} \cap \{\tilde{n}_1,\dots,\tilde{n}_{l_2}\}$ contains at least $\gamma N$ elements.
By \eqref{rectimes} and \eqref{aiutami}, every element of this set is a $(e^{\frac{\lambda}{8}},\delta)$-hyperbolic time, which concludes the proof.

\end{TheoremproofC}
\section{Quenched large deviation estimates for hyperbolic times}
Let $\kappa_1,\delta$ be as in Proposition \ref{hyperbolictimes}. Suppose  $\delta_1$ is such that the claim of  Proposition \ref{hyperbolic ball} holds.\footnote{The precise value of $\delta_1$ is determined in section $5$.}
\begin{definicao}\label{intervalHyperbolic}
Given an interval $I$, we say that {\em $n$ is an  $(\omega,I)$-hyperbolic time } if there exists $x \in I$ such that $n$ is a $(\kappa_1,\delta)$-hyperbolic time for $x$ and $B(x,\delta_1 \kappa_1^{-\frac{n}{2}}) \subset I$.
\end{definicao}
Let $E_N(\omega)$ be the set of $x \in \mathbb{S}^1$ for which the assumptions of Theorem \ref{hyperbolictimes} fail, i.e for which we can not guarantee the density of $(\kappa_1,\delta)$-hyperbolic times in $\{1,\dots,N\}$ for any $(\omega,x)$ such that $x \in E_N(\omega)$:
\begin{align*}
E_N(\omega)&:= \{x \in \mathbb{S}^1 \mid S_N(\omega,x)< \lambda N \,\, \text{or}\,\,Z_{N}(\omega,x,\delta) \ge H(\delta)N \} \qquad \forall N \ge 0,
\end{align*}
where $S_N$ and $Z_N$ are as in Theorem
\ref{hyperbolictimes}.
\begin{proposicao}\label{newTail}
Let $I \subset \mathbb{S}^1$ be an interval. Suppose that $N$ is a natural number satisfying $\gamma N > 2\log_{\kappa_1}\left(\frac{4\delta_1}{|I|}\right)$, where $\gamma $  is as in Theorem \ref{hyperbolictimes}. Suppose further that  $\omega \in \Omega$  is such that $\mathrm{Leb}\left\{E_N(\omega)<\frac{|I|}{2}\right\}$. Then, in the set $\{1,\dots,N\}$ there are at least 
$\gamma N- 2\log_{\kappa_1}\left(\frac{4\delta_1}{|I|}\right)$  $(\omega,I)$-hyperbolic times.
\end{proposicao}
\begin{proof}
Let $y$ be the center of $I$. 
If $\mathrm{Leb}\{E_N(\omega)<\frac{|I|}{2}\}$, then by Theorem \ref{hyperbolictimes} there exists at least one $x \in B_{\frac{|I|}{4}}(y)$ which has  no fewer than $\gamma N$  $(\kappa_1,\delta)$-hyperbolic times in $\{1,\dots,N\}$. Suppose that $n$ is one of such $(\kappa_1,\delta)$-hyperbolic times, and that 
\begin{align*}
n > 2\log_{\kappa_1}\left(\frac{4\delta_1}{|I|}\right).
\end{align*}
Then $B_{\delta_1\kappa_1^{-\frac{n}{2}}}(x) \subset I$, so that $n$ is an $(\omega,I)$-hyperbolic time. This proves the claim.
\end{proof}

The following Proposition is used extensively in the proof of Theorem \ref{tea}, and establishes a bound for the tail of the events  $\left\{\mathrm{Leb}(E_N(\omega))> \frac{|I|}{2} \right\}$.
\begin{proposicao}\label{dev}
 For any interval $I$, there exists  $C_5 \ge 1$, depending on $|I|$,  such that
\begin{equation}\label{final deviation}
\mathbb{P}\left\{\mathrm{Leb}(E_N(\omega))> \frac{|I|}{2} \right\}  < C_5 \left(\frac{1}{N}\right)^{3+\frac{\alpha}{2}} \qquad \forall N \ge 0,
\end{equation}
where $\alpha$ is as in Hypotheses \ref{(L)}-\ref{(M)}.
\end{proposicao}

\begin{proof}
The compiutations are  similar to \cite[Corollary 7.6]{ASENS_2002_4_35_1_77_0} and \cite[Lemma 7.7]{ASENS_2002_4_35_1_77_0} but the exponential tail is replaced by a  polynomial tail.

We first show that there exists $C_3 \ge 1$ such that
\begin{equation}\label{switch}
\mathbb{P}\left\{ \mathrm{Leb}(E_N(\omega))> \left(\frac{1}{N}\right)^{\frac{\alpha}{2}} \right\} \le C_3\left(\frac{1}{N}\right)^{4+\frac{\alpha}{2}} \qquad \forall N \ge 1.
\end{equation}
Indeed, suppose the converse is true and there exists $N_0\ge 1$ such that 
\begin{equation}\label{helpduality}
\mathbb{P}\left\{ \mathrm{Leb}(E_N(\omega))> \left(\frac{1}{N_0}\right)^{\frac{\alpha}{2}} \right\} \ge C_3\left(\frac{1}{N_0}\right)^{4+\frac{\alpha}{2}}.
\end{equation}
For $x \in \mathbb{S}^1$, consider the set of noise realizations
\begin{align*}
E_N(x):=\{ S_N(\omega,x)< \lambda N \,\, \text{or}\,\,Z_{N}(\omega,x,\delta) \ge H(\delta)N\}.
\end{align*}
By Hypotheses \ref{(L)} and \ref{(M)}, there exists $C_3>1,\alpha>0$ such that
\begin{align*}
\mathbb{P}(E_N(x)) \le \frac{C_3}{N^{4+ \alpha}} \qquad \forall N \ge 0, x \in \mathbb{S}^1.
\end{align*}
Thus, by Fubini's Theorem,
\begin{equation}\label{duality}
\int_{\Omega} \mathrm{Leb}(E_{N_0}(\omega)) d\mathbb{P}(\omega)=\int_{\mathbb{S}^1}\mathbb{P}(E_{N_0}(x))d(x) \le C_3\left(\frac{1}{N_0}\right)^{4+\alpha}.
\end{equation}
However, by \eqref{helpduality}
\begin{align*}
\int_{\Omega} \mathrm{Leb}(E_{N_0}(\omega)) d\mathbb{P}(\omega) > \left(\frac{1}{N_0}\right)^{\frac{\alpha}{2}} \mathbb{P}\left\{  \mathrm{Leb}(E_N(\omega)) > \left(\frac{1}{N_0}\right)^{\frac{\alpha}{2}}\right\}> C_3\left(\frac{1}{N_0}\right)^{4+\alpha},
\end{align*}
which contradicts \eqref{duality}, and this implies \eqref{switch}.


Now, consider
\begin{align*}
\tilde{E}_k:= \left\{ \exists i \ge k\,\text{such that}\, \mathrm{Leb}(E_i(\omega)) > \left(\frac{1}{i}\right)^{\frac{\alpha}{2}}\right\}.
\end{align*}
By \eqref{switch}, we have
\begin{align*}
\mathbb{P}(\tilde{E}_k) \le C_3\sum_{i=k}^{\infty} \left(\frac{1}{i}\right)^{4+\frac{\alpha}{2}}
\le C_4\left(\frac{1}{k}\right)^{3+\frac{\alpha}{2}},
\end{align*}
for some constant $C_4 \ge 1$.
Since  $\lim_{n \to \infty} \mathbb{P}(\tilde{E}_n)=0$, the random time
\begin{align*}
N_1(\omega):= \min \{n \ge 1 \mid \omega \notin \tilde{E}_n \}
\end{align*}
is finite almost surely. By the definition of $\tilde{E}_{N_1(\omega)}$, for all  $N \ge N_1(\omega)$
\begin{equation}\label{tr1}
\mathrm{Leb}(\tilde{E}_N(\omega)) \le \left(\frac{1}{N}\right)^{\frac{\alpha}{2}}.
\end{equation}
Furthermore,
\begin{equation}\label{tr2}
\mathbb{P}\{ N_1(\omega) > k\} \le C_4\left(\frac{1}{k}\right)^{3+\frac{\alpha}{2}}.
\end{equation}

Finally, for $k_0$ sufficiently large so that $\left(\frac{1}{k_0}\right)^{\frac{\alpha}{2}}<\frac{|I|}{2}$, by \eqref{tr1}, if $\mathrm{Leb}(E_N(\omega))> \frac{|I|}{2} $, then $N_1(\omega)>N-k_0 $
and \eqref{final deviation} follows from \eqref{tr2}.
\end{proof}
\section{Random Horseshoe}
In this section we prove Theorem \ref{tea}. The main ingredients in the proof of this theorem are Lemma \ref{return times}, Proposition \ref{almostDone} and Proposition \ref{PosD}. 

Let $\kappa_0$ be as in Hypothesis \ref{(O)} and $\kappa_1$ as in Theorem \ref{hyperbolictimes}, take $\kappa:= \min\{\kappa_0,\kappa_1\}$. Let $I$ be any given interval and consider the following random time
\begin{equation}\label{randomT}
m(\omega,I):= \min\{n \ge 0 \mid I \,\,\text{is an}\,\,(\omega,n,\kappa)\text{-full expanding  branch}\}. 
\end{equation} 
In the following Lemma we show that the random time $m$ is almost surely finite with at least a power-law tail.
\begin{lema}\label{return times}
For every  interval $I$, there exists a constant $C_7 \ge 1 $ such that
\begin{equation}\label{return times eq}
\mathbb{P}\{ m(\omega,I)> l \} \le C_7 \left(\frac{1}{l}\right)^{3+\frac{\alpha}{4}} \qquad \forall l \ge 0.
\end{equation}
In particular, $m(\omega,I)$ is almost surely finite. 
\end{lema}
We defer the proof of this Lemma to the end of the section, as it is involved.

Define recursively the random sequence
\begin{equation}\label{defrecurs}
n_{k}(\omega):= \max\{m(\theta^{n_{k-1}(\omega)}(\omega),I_0),m(\theta^{n_{k-1}(\omega)}(\omega),I_1) \}.
\end{equation}
By Lemma \ref{return times},  $n_k$ is almost surely finite for all $k \ge 0$. In the following Proposition, we show how to use the definition of $\{n_k(\omega)\}_{k \ge 0}$ to construct symbolic dynamics between $I_0$ and $I_1$.
\begin{proposicao}\label{almostDone}
Given two disjoint  intervals $I_0,I_1 \subset \mathbb{S}^1$, there exists a full measure set $\Hat{\Omega}$, such that for all $\omega \in \Hat{\Omega}$, $i,j \in \{0,1\}$ and $k \ge 0$, there exists a set  $J^{i,j}_k(\omega) \subset I_i$ such that
\begin{align*}
 g^{n_{k+1}(\omega)-n_k(\omega)}_{\theta^{n_k(\omega)}(\omega)}(J^{i,j}_k(\omega))&= I_j,\\
 \left|dg^{n_{k+1}(\omega)-n_k(\omega)}_{\theta^{n^k(\omega)}(\omega)}|_{J^{i,j}_k(\omega)}\right| &> \kappa > 1.
\end{align*}
\end{proposicao}
\begin{proof}
Let $\Hat{\Omega}$ be such that $n_k$ is finite for  all $\omega \in \Hat{\Omega}$ and all $k \ge 0$. Consider the case $k= 0$. Without loss of generality let  $m(\omega,I_0) > m(\omega,I_1)$.
Let $E_{\omega} := \{x \in \mathbb{S}^1 \mid |dg_{\omega}|> 1\}$.
By Hypothesis \ref{(P)}, up to restricting to a full measure subset of noise realizations, we can assume that $g_{\theta^s(\omega)}(E_{\theta^s(\omega)})= \mathbb{S}^1$ for all $\omega \in \Hat{\Omega}$ and all $s \ge 0$. As a consequence, if $\omega \in \Hat{\Omega}$ there are sets $E'_{{\omega},j} \subset E_{\theta^{m(\omega,I_1)}(\omega)}$, $j \in \{0,1\}$, such that $g^{q}_{\theta^{m(\omega,I_1)}(\omega)} \subset E_{\theta^{m(\omega,I_1)+i}(\omega)}$ for every $q \le m(\omega,I_0)-m(\omega,I_1)$ and
\begin{equation}\label{circlecover}
g^{m(\omega,I_0)-m(\omega,I_1)}_{\theta^{m(\omega,I_1)}(\omega)}(E'_{{\omega},j})= I_j.
\end{equation}
By definition of $m(\omega,I)$, we can define the sets $\{J^{i,j}_0(\omega)\}_{i,j \in \{0,1\}}$ such that $J^{i,j}_0(\omega) \subset I_i$ and 
\begin{align*}
 g^{m(\omega,I_1)}_{\omega}(J^{1,j}_0(\omega))&= E'_{{\omega},j},\\ 
  g^{m(\omega,I_0)}_{\omega}(J^{0,j}_0(\omega))&= I_j, \\
 \left|dg^{m(\omega,I_i)}_{\omega}\vert_{J^{i,j}_{0}(\omega)}\right| &> \kappa. 
\end{align*}
By \eqref{circlecover}
\begin{align*}
g^{m(\omega,I_0)}_{\omega}(J^{1,j}_0(\omega))&= g^{m(\omega,I_0)-m(\omega,I_1)}_{\theta^{m(\omega,I_1)}(\omega)}\left(g^{m(\omega,I_1)}_{\omega}(J^{1,j}_0(\omega))\right)\\
&= g^{m(\omega,I_1)-m(\omega,I_0)}_{\theta^{m(\omega,I_1)}(\omega)}(E'_{\omega,j})= I_j.
\end{align*}
This concludes the proof for $k=0$ because we assumed that $n_0(\omega)=m(\omega,I_0)$. For the case $k \ge 1$ the proof is the same up to replacing $\omega$ with $\theta^{n_{k-1}(\omega)}(\omega)$.
\end{proof}
In the following Proposition, we establish positive density for the sequence $\{n_k(\omega)\}_{k \ge 0}$.
\begin{proposicao}\label{PosD}
 For every pair of disjoint intervals $I_0,I_1 \subset \mathbb{S}^1$, there exists a full measure set $\overline{\Omega} \subset \Hat{\Omega}$ such that
\begin{align*}
\lim_{k \to \infty}\frac{n_k(\omega)}{k}= E[n_0(\omega)] \qquad  \forall \omega \in \overline{\Omega}.
\end{align*}
\end{proposicao}
\begin{proof}
Observe that 
\begin{align*}
n_k= n_0+ (n_1-n_0)+ \dots + (n_{k}-n_{k-1}).
\end{align*}
The random variables $\{n_{j+1}-n_j\}_{j \ge 0}$ and $n_0$ are independent and identically distributed. Furthermore, by  \eqref{return times eq}
\begin{align*}
\mathbb{E}[(n_0)^2] &\le \sum_{l \ge 0} l^2 \mathbb{P}(n_0(\omega) > l-1) \\&\le C_6 \sum_{l \ge 0}l^2 \left(\frac{1}{l}\right)^{3+ \frac{\alpha}{4}}\\
&\le C_6\sum_l \frac{1}{l^{1+\frac{\alpha}{4}}} < \infty.
\end{align*}
By the law of large numbers, we conclude that
\begin{align*}
 \lim_{k \to \infty}\frac{n_k(\omega)}{k}= \mathbb{E}[n_0(\omega)], \qquad \mathbb{P} \,\,\text{a.s.}
\end{align*}
\end{proof}
\begin{proof}[Proof of Theorem \ref{tea}]
Fix two disjoint intervals $I_0,I_1 \subset \mathbb{S}^1$ and consider the random sequence $\{n_k(\omega)\}_{k \ge 0}$. By Proposition \ref{almostDone} and \ref{PosD} there is a full $\mathbb{P}$-measure set $\overline{\Omega}(I_0,I_1)$ such that $\{n_k(\omega)\}_{k \ge 0}$ satisfies \eqref{positiveDensity}, \eqref{e2} and \eqref{tea} for all $\omega \in \overline{\Omega}(I_0,I_1)$. By Definition \ref{Horseshoe}, the pair $(I_0,I_1)$ is an $(\omega,\kappa)$-horseshoe for all $\omega \in \overline{\Omega}(I_0,I_1)$.  
This implies the theorem with the full measure set
\begin{align*}
\tilde{\Omega}:= \bigcap _{I_0,I_1 \,\text{dyadic intervals}} \overline{\Omega}(I_0,I_1).
\end{align*}
\end{proof}
The rest of this section is dedicated to the proof of Lemma \ref{return times}.
\begin{proof}[Proof of Lemma \ref{return times}]
Fix a interval $I \subset \mathbb{S}^1$. Let $\delta_1>0$ be small enough such that there exists $K$ and $\iota$ as in Hypothesis \ref{(O)}. Define the sequence of $K$-sparse hyperbolic times
\begin{align*}
\tau_1(\omega)&:= \min \{ n \ge 1 \mid \text{$n$ is a $(\omega,I)$-hyperbolic time } \},\\
\tau_{k+1}(\omega)&:=\min \{n > \tau_k(\omega)+K \mid \text{$n$ is a $(\omega,I)$-hyperbolic time} \},
\end{align*}
where the $(\omega,I)$-hyperbolic times are as in Definition \ref{intervalHyperbolic}.
Let $k_l(\omega):= \max \{ i \mid \tau_i(\omega) \le l\}$, i.e. $k_l(\omega)$ is the number of $K$-sparse hyperbolic times that lie in $\{1,\dots,l\}$. Before proceeding with  the proof, in order not to interrupt the flow of ideas, we establish an estimate which will be used in the conclusion of this section. 
\begin{proposicao}\label{kestimate}
Let $s$ be a natural number such that  $s(K+1)   \le  \max\left\{0,\gamma l- 2\log_{\kappa}\left(\frac{4\delta_1}{|I|}\right)\right\}$, where $\gamma \in (0,1)$ is as in Theorem \ref{hyperbolictimes}. Then 
\begin{align*}
\mathbb{P}(k_l(\omega) = s) \le C_6 \left(\frac{1}{l}\right)^{3+\frac{\alpha}{2}}.
\end{align*}
 
\end{proposicao}
\begin{proof}
Let $t_l(\omega)$  be the number of $(\omega,I)$-hyperbolic times less or equal to $l$. It is readily verified that
\begin{equation}\label{tijmans}
k_l(\omega) \ge \frac{1}{K+1}t_l(\omega), \qquad \forall l \ge 0.
\end{equation}
By Proposition \ref{newTail}, if $s(K+1) \le  \max \left\{ \gamma l - 2\log_{\kappa}\left(\frac{4\delta_1}{|I|}\right),0 \right\}$, then either $ l \le \frac{2\log_{\kappa}\left(\frac{4\delta_1}{|I|}\right)}{\gamma}$ or
\begin{align*}
\left\{k_l(\omega) \le  s \right\} \subset  \left\{ \mathrm{Leb}(E_l(\omega))> \frac{|I|}{2} \right\},
\end{align*}
and the result follows immediately from Proposition \ref{dev} with $l=N$.
\end{proof}

By definition of the $(\omega,I)$-hyperbolic times, the image  $f^{\tau_j}_{\omega}(I) $ contains a ball of radius $\delta_1$ for all $j \in \mathbb{N}$. Let $H(\omega,j)$ be such ball.
Define the set of noise realizations
\begin{align*}
\mathcal{S}\left(K,j,\kappa\right):=  \{ H(\omega,j) \,\text{is not } (\theta^{\tau_j}(\omega),K,\kappa)\text{- full branch expanding}\},\qquad j \ge 0.
\end{align*}
Observe that if $m(\omega,I)>l$, then $I$ is not $(\omega,l,\kappa)$-full branch expanding. This implies  that $H(\omega,j)$ is  not $(\tau_j(\omega),K,\kappa)$-full branch expanding, for all $j = 1,\dots, k_l(\omega)-1$.
Therefore, we can write
\begin{align*}
\mathbb{P}\left(m(\omega,I) > l\right)&=\sum_{s=0}^l\mathbb{P}\left(m(\omega,I) > l\,\,\text{and}\,\,k_l(\omega)=s\right)\\ &\le \sum_{s=0}^l\mathbb{P}\left(\bigcap_{j=1}^{k_l(\omega)-1} \mathcal{S}\left(K,j ,\kappa\right) \,\,\text{and}\,\, k_l(\omega)=s \right)
\\&= \sum_{s=0}^l   \mathbb{P}\left(\bigcap_{j=1}^{s-1} \mathcal{S}\left(K,j,\kappa\right) \,\,\text{and}\,\, k_l(\omega)=s \right).
\end{align*}
Therefore, by  H\"older inequality, we have 
\begin{equation}\label{jugete}
\mathbb{P}\left(m(\omega,I)>l \right) \le  \sum_{s=0}^l \left({\mathbb{P} \left(\bigcap_{j=1}^{s-1} \mathcal{S}\left(K,j,\kappa\right)\right)}\right)^{\frac{\alpha}{12+2\alpha}}\left({\mathbb{P}(k_l=s)}\right)^{\frac{12+\alpha }{12+2\alpha}},
\end{equation}
where $\alpha$ is as in Proposition \ref{dev}.

Let $\iota$ be as in Hypothesis \ref{(P)}. We claim that
\begin{equation}\label{inequality}
\mathbb{P}\left(\bigcap_{j=1}^{s-1} \mathcal{S}\left(K,j,\kappa\right)\right) \le (1-\iota)^{s-1}.
\end{equation}
Let $V_i= \bigcap_{j=1}^i\mathcal{S}\left(K,j,\kappa\right)$, then 

\begin{equation}\label{mistery}
\mathbb{P}\left(V_{s-1}\right)=\mathbb{E}\left[\mathbb{P}\left(\mathcal{S}\left(K,s-1,\kappa\right)|\mathcal{F}_{\tau_{s-1}}\right)(\omega)\cdot1_{V_{s-2}}(\omega)\right],
\end{equation}
where $\mathcal{F}_{\tau_{s-1}}$ denotes the past $\sigma$-algebra up to time $\tau_{s-1}$.
 Since the length of $H(\omega,s-1)$ is equal to $2\delta_1 \in (\frac{\eta}{4},4\eta)$, then by Hypothesis \ref{(O)},
\begin{align*}
 \mathbb{P}\left(\mathcal{S}\left(K,s-1,\kappa\right)|\mathcal{F}_{\tau_{s-1}}\right)(\omega) \le 1-\iota.
 \end{align*}
 Therefore, by \eqref{mistery}
\begin{align*}
\mathbb{P}\left(\bigcap_{j=1}^{s-1}\mathcal{S}\left(K,\tau_j,H(\omega,j) ,\kappa\right)\right) \le  (1-\iota)\mathbb{P}\left(\bigcap_{j=1}^{s-2}\mathcal{S}\left(K,\tau_j,H(\omega,j) ,\kappa\right)\right),
\end{align*}
and inequality \eqref{inequality} follows by induction.
It follows from equations \eqref{jugete}-\eqref{inequality} and Proposition \ref{kestimate} that 
\begin{align*}
\mathbb{P}(m(\omega,I)>l) &\le (1-\iota)^{\frac{\alpha}{12+2\alpha}}\sum_{s = 0}^l (1-\iota)^{s\frac{\alpha}{12+2\alpha}} {\mathbb{P}(k_l(\omega)=s)}^{\frac{12+\alpha }{12+2\alpha}}\\
&\le  (1-\iota)^{\frac{\alpha }{12+2\alpha}}C_6\left(\frac{1}{l}\right)^{\left(3+\frac{\alpha}{2}\right)(\frac{12+\alpha}{12+2\alpha})}\sum_{s(K+1) \le  \max\left\{0,\gamma l- 2\log_{\kappa}\left(\frac{4\delta_1}{|I|}\right)\right\}}  (1-\iota)^{s\frac{\alpha}{12+2\alpha}}\\
&+(1-\iota)^{\frac{\alpha}{12+2\alpha}}\sum_{s(K+1) > \max\left\{0,\gamma l- 2\log_{\kappa}\left(\frac{4\delta_1}{|I|}\right)\right\}}  (1-\iota)^{s\frac{\alpha}{12+2\alpha}}  \\
&\le C_6\left(\frac{1}{l}\right)^{3+\frac{\alpha}{4}}(1-\iota)^{\frac{\alpha}{12+2\alpha}}\left(\sum_{s \ge 0} (1-\iota)^{s\frac{\alpha}{12+2\alpha}}\right)\\ &+ \left((1-\iota)^{\frac{-\alpha\log\left(\frac{4\delta_1}{|I|}\right)}{12+2\alpha}}\right)\left((1-\iota)^{\frac{\alpha}{12+2\alpha} \frac{\gamma}{K+1}}\right)^l\left(\sum_{s \ge 0} (1-\iota)^{s\frac{\alpha}{12+2\alpha}}\right)  \\
&\le C_7\left(\frac{1}{l}\right)^{3+\frac{\alpha}{4}},
\end{align*}
for some constant $C_7 \ge 1$, which concludes the proof of Lemma \ref{return times}.
\end{proof}
\section{A random horseshoe for predominantly expanding RDS: Lyapunov exponent and invariant measure}
Section $6$ and $7$ are  devoted to the proof of Theorem \ref{teb}. In particular, we need to verify  that any $(\sigma,R)$-predominantly expanding RDS (see Definition \ref{admissibleRDS})  satisfies Hypotheses \ref{(A)}-\ref{(P)}. In Section $6$ we check that any $(\sigma,R)$-predominantly expanding RDS satisfies Hypotheses \ref{(A)}-\ref{(U)} and \ref{(O)}-\ref{(P)}, whilst in Section $7$ we focus on the remaining Hypotheses \ref{(L)}-\ref{(M)}. 

From now on, using the same notation as in Section $2.3$, we set  $\Omega:= [-\sigma,\sigma]^{\mathbb{N}}$, $\mathbb{P}:= \left(\frac{\mathrm{Leb}\vert_{[-\sigma,\sigma]}}{2\sigma }\right)^{\otimes \mathbb{N}}$  given the noise amplitude $\sigma>0$ and we let $\mathcal{F}$ be the  product Borel $\sigma$-algebra on $\Omega$ obtained via Kolmogorov extension theorem.

We remind-  see Section $2.3$  for further details-  that a $(\sigma,R)$-predominantly expanding RDS $f^n_{\omega}$ is generated by a family of iterations $\{f_{\omega}\}_{\omega \in \Omega}$ that has  a regular random singular set $\{\mathcal{SC}_{\omega}\}_{\omega \in \Omega}$ (see Definition \ref{cset}) and satisfies
\begin{align*}
f_{\omega}(x):= f(x + \omega_0), 
\end{align*}
for a given function $f$ that satisfies the requirements in Definition \ref{admissibleRDS}. Due to the behaviour of the first and the second derivative of $f$ around the singular set $\mathcal{S}$, item $(ii)$ of Definition \ref{admissibleRDS}, one can prove that  $\mathcal{SC}$ is regular (see Definition \ref{cset}). Recall that $\mathcal{SC}$  is the union of the non-differentiability set $\mathcal{S}$ and of the critical set $\mathcal{C}$. Furthermore, recall that $G= \{|df|>R \}$ is the super expanding set, and there are $G_1,\dots,G_k$ connected components of $G$ that satisfy $|G_i| \ge \frac{1}{R}$. 

The following proposition shows a crucial property of $(\sigma,R)$-predominantly expanding maps, which is frequently used in this section, i.e. that there exists a $\Delta$ in the super-expanding region such that $|\Delta|\le \frac{1}{R}$, $f(\Delta)= \mathbb{S}^1$ and $\Delta$ satisfies an accessibility condition (see equation \eqref{elefante}).
\begin{proposicao}\label{newresult}
Let $f^n_{\omega}$ be $(\sigma,R)$-predominantly expanding. Let $G_1$ be one of the connected components of $G$ in item $(iii)$ of Definition \ref{admissibleRDS} and let $c$ be its middle point. Let 
\begin{equation}\label{reference}
\Delta:= \left[c-\frac{1}{2R},c+\frac{1}{2R}\right].
\end{equation}
Then, $\Delta \subset G_1$ and $f(\Delta)= \mathbb{S}^1$.

Furthermore, given $\gamma>0$, consider the sets
\begin{equation}\label{gamma}
\begin{array}{ll}
\Delta(\gamma) &:= \{ x \in \Delta \mid d\left(x,\mathbb{S}^1 \setminus G\right)>\gamma \}, \\
E(\gamma) &:= \{ x \in E \mid d\left(x,\mathbb{S}^1 \setminus E\right)>\gamma \}.
\end{array}
\end{equation}
Then, there exists a $\gamma>0$, $M \ge 0$ and $q \in (0,1)$ such that
\begin{equation}\label{elefante}
\mathbb{P}\left( x+\omega_0 \in E(\gamma)\, \,\,\text{and}\,\, f_{\omega}(x) \in \Delta(\gamma)\right) > q,
\end{equation}
for all $x \in \mathbb{S}^1$.
\end{proposicao}\label{youliar}
\begin{proof}
The first assertion of this proposition follows from 
the fact that  any interval inside $G_1$ of size at least $\frac{1}{R}$ is mapped into the full circle due to Lagrange theorem. For the proof of \eqref{elefante}, observe that, since  $\mathrm{Leb}\left(\mathbb{S}^1 \setminus G\right)<D(R)$ by \eqref{dierre}, then for any $x \in \mathbb{S}^1$ there exists $G_{i}$ such that $d(x,G_i)< \frac{D(R)}{2}$. Using\eqref{noiseRange} and the fact that $f(G_i)= \mathbb{S}^1$, choosing $\gamma$  small enough such that
\begin{align*}
\left(G_i \cap \Delta(\gamma)\right) \cap f^{-1}(\Delta(\gamma)) \neq \emptyset\,\, \forall i =1,\dots, k, 
\end{align*}
we conclude that
\begin{align*}
\inf_{x \in \mathbb{S}^1}\mathbb{P} \left(x+\omega_0 \in f^{-1}(\Delta(\gamma))\cap G_i \right)>0,
\end{align*}
from which \eqref{elefante} follows immediately.
\end{proof}
\subsection{Regularity of the singular set and uniqueness of ergodic stationary measure}

The main result of this subsection is the following proposition:
\begin{proposicao}\label{SoyAquel}
Let $f^n_{\omega}$ be any $(\sigma,R)$-predominantly expanding random dynamical system. Then \ref{(A)} and \ref{(U)} are satisfied. In particular, the RDS $f^n_{\omega}$ has an absolutely continuous unique ergodic stationary measure $\mu$ such that  $\log|df_{\omega}(x)| \in L^1(\mathbb{P} \otimes \mu)$. 
\end{proposicao}
\begin{proof}
The proof that \ref{(A)} is verified is an easy consequence of item $(i)$ in Definition \ref{admissibleRDS} and is left to the reader.

 To prove existence and uniqueness of the ergodic measure, due to compactness of $\mathbb{S}^1$, it is sufficient to prove that  any initial condition $x \in \mathbb{S}^1$ can reach any open set of the circle  with positive probability ( see \cite[chapter 1]{kifer2012ergodic}). For a $(\sigma,R)$-predominantly expanding RDS $f^n_{\omega}$ this fact  is established in the  following Lemma,  which uses the existence of  $\Delta$ as in \eqref{newresult} and the sufficiently large noise amplitude. 
\begin{lema}\label{amantes}
For any $x \in \mathbb{S}^1$ and  for any open set $A \subset \mathbb{S}^1$, $\mathbb{P}\left\{f_{\omega}^2(x) \in A\right\}>0$.
\end{lema}
\begin{proof} 
Fix $x \in \mathbb{S}^1$ and take an open set $A \subset \mathbb{S}^1$. By \eqref{elefante}
\begin{equation}\label{fullb}
\mathbb{P}\left\{f_{\omega}(x) \in \Delta \right\}>0.
\end{equation}
Since $f(\Delta) = \mathbb{S}^1$, $f^{-1}(A) \cap \Delta \neq \emptyset$. Then,
\begin{align*}
\mathbb{P}\left\{ f^{2}_{\omega}(x) \in A\right\} &\ge  \mathbb{P}\left\{ f_{\omega}(x) +\omega_{1} \in \Delta \cap f^{-1}(A)\right\}  \\
&\ge  \mathbb{P}\left\{f_{\omega}(x) \in \Delta \right\}\inf_{y \in \Delta}\mathbb{P}\left\{ y+\omega_0  \in f^{-1}(A) \cap \Delta  \right\}.
\end{align*}
In the last line, the first term is greater than zero due to \eqref{fullb}, the second term is positive due to the fact that $\sigma > \frac{1}{R}= |\Delta|$. This concludes the proof.
\end{proof}

As a consequence of Lemma \ref{amantes}, the RDS $f^n_{\omega}$ has an unique ergodic measure $\mu$. The absolute continuity of $\mu$ follows from the fact that the transition probabilities $A \in \mathcal{B}(\mathbb{S}^1) \mapsto \mathbb{P}(f_{\omega}(x) \in A)$  are absolutely continuous with respect to Lebesgue for all $x \in \mathbb{S}^1$ \cite[chapter 1]{kifer2012ergodic}.

It remains to prove that $\log|df_{\omega}(x)|$ is integrable with respect to $\mathbb{P}\otimes \mu$. 
Observe that $df(x)$ is uniformly bounded from above and is bounded away from $0$ everywhere but in a small neighbourhood of the critical set $\mathcal{C}$, in which $df(x)$ approaches $0$. As a consequence, we only need to check that the integral of  $\log|df_{\omega}(x)|$ with respect to $\mathbb{P}\otimes \mu$ is bounded from below.
Due to the regularity of $\mathcal{SC}$, there exists a constant $B \ge 1$ such that, for all $x$ in a small neighbourhood $\mathcal{N}(\mathcal{C})$ of $\mathcal{C}$
\begin{align*}
|df(x)| \ge B \mathrm{dist}(x,\mathcal{C}).
\end{align*}
Using this fact, we see that, there exists a constant $D_1 >0$ such that for all $x \in \mathbb{S}^1$
\begin{align*}
\int_{\Omega}\log|df(x+\omega_0)|d\mathbb{P}(\omega_0) &\ge \int_{B(x,\sigma)} \log\left|df(z)\right| dz\\
&\ge \int_{B(x,\sigma) \setminus \mathcal{N}(\mathcal{C})}\inf_{ x \in \mathbb{S}^1\setminus \mathcal{N}({\mathcal{C})}}\log|df(x)| 
\\ &+ \int_{\mathcal{N}(\mathcal{C})}
\log B \mathrm{dist}(z,\mathcal{C})dz \\
&\ge  D_1 + \int_{\mathcal{N}(\mathcal{C})}
\log\mathrm{dist}(z,\mathcal{C})dz.
\end{align*}
Therefore
\begin{align*}
\int_{\mathbb{S}^1\times \Omega}\log|df(x+\omega_0)|d\mu(x)d\mathbb{P}(\omega_0) &= \int_{\mathbb{S}^1}\left[\int_{\Omega} \log\left|df(x+\omega_0)\right|d\mathbb{P}(\omega)\right] dx
\\
&\ge D_1+ \int_{\mathcal{N}(\mathcal{C})} \log\mathrm{dist}(x,\mathcal{C})dx>-\infty,
\end{align*}
which concludes the proof.
\end{proof}
\subsection{Full branch expanding property}
We have proved that Hypotheses \ref{(A)} and \ref{(U)} are fulfilled. Note also that Hypothesis \ref{(P)} is trivially satisfied. Indeed let $G_{\omega}$ be the super expanding region of $f_{\omega}$ and consider also the sets  $\Delta_{\omega}:= \Delta-\omega_0$, for $\omega \in \Omega$. Then $\Delta_{\omega} \subset G_{\omega}$ and $f_{\omega}(\Delta_{\omega})= \mathbb{S}^1$, for all $\omega \in \Omega$, due to Definition \ref{admissibleRDS}. Now we prove that \ref{(O)} is also fulfilled. \\
Recall that (see Definition \ref{fullBranchExpanding}) an interval $I$ is $(\omega,\kappa_1,K)$-full branch expanding if 
there exist an interval $J \subset I$ and an integer $i \le K$ such that $f^i_{\omega}(J)= \mathbb{S}^1$ and $\left|df^i_{\omega}\vert_{J} \right|>\kappa_1$.

The main result of this subsection is the following 
\begin{proposicao}\label{conmigoestabachata}
Let $f^n_{\omega}$ be a $(\sigma,R)$-predominantly expanding  RDS and let $G$ be the super-expanding region of $f$. Let $\Delta \subset G$ such that $f(\Delta)=\mathbb{S}^1$. Then, Hypothesis \ref{(O)} is  fulfilled, namely
for sufficiently small $\eta >0$ there exists numbers  $K \in \mathbb{N}$, and $\kappa_1 \ge 1$ such that there is a positive probability, bounded away from zero, that any interval whose diameter lies in $(\frac{\eta}{4},4\eta)$ is $(\omega,\kappa_1,K)$-full branch expanding. In particular, \eqref{ctrol} follows.
\end{proposicao}
The proof is divided into three steps.
\begin{itemize}
    \item In Lemma \ref{sincensura}, we use \ref{(P)} to prove that for all $k \in \mathbb{N}$ and all $\omega \in \Omega$ there exists a random  interval $J=J(\omega_0,\dots,\omega_{k-1})$ such that $J$ is $(\omega,k,R^k)$-full branch expanding and $|J| \le R^{-k}$. 
    \item In Lemma \ref{hacecalor}, we deduce from Lemma \ref{sincensura} uniform bounds on the probability that, given any $\eta_1>0$,  there exists  $K \in \mathbb{N}$  such that any interval  $I$ contained  in the set $\Delta$ introduced in \eqref{newresult}, with $|I| \ge \eta_1$, is $(\omega,K,R^{K})$-full branch expanding; 
    \item We combine Lemma \ref{sincensura} and Lemma \ref{hacecalor}, together with Proposition \ref{reference}, to  conclude the proof.
\end{itemize}
\begin{lema}\label{sincensura}
 For every $k \in \mathbb{N}$, there exist 
 $\gamma_k \in (0,1)$ and a set of noise realizations $H_k$ of the form
 \begin{equation}\label{fixing}
 H_k:=[-\sigma,\sigma]\times H'(\omega_1,\dots,\omega_{k-1})
 \end{equation}
 satisfying $\mathbb{P}(H_k)>\gamma_k$ and the property that for  every $\omega \in H_k$ there exists a random  interval  $J_k=J_k(\omega_0,\dots,\omega_{k-1})$  such that 
$|J_k| \le R^{-k}$, $f^i_{\omega}(J_k) \subset \Delta_{\theta^{i}(\omega)}$ for all $0= 1,\dots,k-1$ and $f^k_{\omega}(J_k)=\mathbb{S}^1$, where $\Delta_{\omega}:= \Delta-\omega_0$, and $\Delta$ is as in \eqref{newresult}.  In particular, $J_k$ is  $(\omega,k,R^k)$-full branch expanding.
\end{lema}
\begin{proof}
We first state, without proof, the following useful fact.
\begin{lema}\label{useful}
Let $Q \subset \mathbb{S}^1$ be a interval such that $|Q|\le \frac{1}{R}$. Then
\begin{equation}\label{cutting bad parts}
\mathbb{P}\left\{ f^{-1}(Q+\omega_0)\cap \Delta \,\,\text{is connected} \right\}\ge \frac{1}{4R\sigma}.
\end{equation} 
\end{lema}
We  can then proceed with the proof of Lemma \ref{sincensura}. Let $\bar{f}$ be the lifting of $f$ and $\bar{\Delta}$ be a lifting of $\Delta$ such that $\bar{\Delta} \subset [0,1)$. Let us write $\bar{\Delta}:=[a,b]$, with $0\le a < b \le 1$ and $b-a\le \frac{1}{R}$. Without loss of generalities, let us assume that $\bar{f}(a)=0$ and $\bar{f}(b)=1$.  Fix $k \in \mathbb{N}$ and define $J_1 := \Delta_{\theta^{k-1}(\omega)}$, so that $f(J_1+\omega_{k-1})= \mathbb{S}^1$. Note that $J_1$ is connected by definition.
Now take $J_2 \subset \Delta_{\theta^{k-2}(\omega)}$ such that $f( J_2+\omega_{k-2})=J_1$. It follows that 
\begin{align*}
J_2:= f^{-1}(J_1+\omega_{k-1})\cap \Delta-\omega_{k-2},
\end{align*}
 is connected with probability greater than $\frac{1}{4R\sigma}$ by Lemma \ref{useful}. 
Proceeding inductively, we construct a sequence of random intervals $\{J_i\}_{i=1,\dots,k}$ such that $J_i \subset \Delta_{\theta^{k-i}(\omega)}$, $f(J_i+\omega_{k-i})=J_{i-1}$ for $i \ge 2 \le k-1$, $f(J_1+\omega_{k-1})= \mathbb{S}^1$ and all intervals $\{J_i\}_{i =1}^n$ are connected as long  as the noise realizations $\{\omega_i\}_{i=1}^{k-1}$ are chosen so that 
\begin{align*}
J_i:= f^{-1}(J_{i-1}+\omega_{k-i}) \cap \Delta 
\end{align*}
is connected for all $i=2,\dots, k-1$. By equation \eqref{cutting bad parts}, this happens with probability larger than $\left(\frac{1}{4R\sigma}\right)^{k-1}$. The other properties required for the random sets $J_{i}$ follow by construction.
\end{proof}
In the following Lemma we show that intervals in $\Delta$ are full branch expanding.
\begin{lema}\label{hacecalor}
For any $\eta_1>0$, there exists $K_1 \in \mathbb{N}$ and $\iota_1 \in (0,1)$ such that, if $|I| \ge \eta_1$ and $I \subset \Delta$, then
\begin{align*}
\mathbb{P}\left\{I \,\,\text{is} \,\, (\omega,K_1,R^{K_1})\text{-full branch expanding}\right\}\ge \iota_1.
\end{align*}
\end{lema}
\begin{proof}
Let $K_1$ be such that $\eta_1 \gg R^{-K_1}$.
 Let $H_k \subset \Omega$ be the set of noise realization defined in \eqref{fixing}. Now consider the set
 \begin{align*}
  \tilde{H}_k:= \theta^{-1}(H_k),
 \end{align*}
 so that for all $\omega \in \tilde{H}_k$ there exists an interval  $J=J(\omega_1,\dots,\omega_{K_1})$
such that $f^{K_1-1}_{\theta(\omega)}(J)= \mathbb{S}^1$, $f^i_{\theta(\omega)}(J) \subset \Delta_{\theta^{i+2}(\omega)}$ and  $|J| \le R^{-K_1+1}$.

Let $\tilde{J}= \tilde{J}(\omega_1,\dots,\omega_{K_1}) \subset \Delta$ be such that $f(\tilde{J})= J$. Since $\tilde{J} \subset \{|df|> R \}$, we have that   $\left|\tilde{J}\right|\le R^{-K_1}<<\eta_1$. Let $\tilde{H}_{k+1}$ denote the set of noise realization such that $\tilde{J}$ is connected. Observe that
\begin{align*}
\left\{ I \,\,\text{is} \,\, (\omega,K_1,R^{K_1})\text{-full branch expanding}\right\} \supset \left\{I+\omega_0 \supset \tilde{J}\,\,\text{and}\,\, \tilde{J}\,\text{is connected} \right\}.
\end{align*}
As a consequence, we see that
\begin{equation}\label{equation}
\begin{array}{l}
\mathbb{P}\left\{I \,\,\text{is} \,\, (\omega,K_1,R^{K_1})\text{-full branch expanding}\right\}\ge \mathbb{P}\left\{ I+ \omega_0 \supset \tilde{J}\,\,\text{and} \,\, \tilde{J} \,\,\text{is connected}\right\}  \\ \displaystyle
 \qquad \ge \int_{\tilde{H}_{k+1}} \left[\int\mathbb{P}\left\{I+\omega_0\ \supset  \tilde{J}   \right\}d\mathbb{P}(\omega_0)\right]d\mathbb{P}(\omega_1,\dots,\omega_{K_1}).
\end{array}
\end{equation}
 Since both $I$ and $\tilde{J}$ are contained in $\Delta$ and  $\sigma> \frac{1}{R}> |\Delta|$ due to the $(\sigma,R)$-predominant expansion, there exists $\iota_1>0$ such that
\begin{align*}
\mathbb{P}\left\{I+\omega_0\,\,\text{contains}\,\, \tilde{J} \right\}>\iota_1,
\end{align*}
for all $\omega_1,\dots,\omega_{K_1} \in \tilde{H}_{k+1}$. By \eqref{equation} and \eqref{fixing} this concludes the proof.




\end{proof}
\begin{TheoremproofE}

Let $C:= \sup_{\omega \in \Omega}||g_{\omega}||_{\mathcal{C}^1\setminus \mathcal{SC}_{\omega}}$. Define $\eta:=\frac{\gamma}{2C}$. Let $I$ be an interval such that $|I| \in (\frac{\eta}{4},4\eta)$. In particular, $I \subset B(x_I,2\eta)$, where $x_I$ is the middle point of $I$. Using \eqref{elefante} with $x=x_I$, we see that 
\begin{align*}
\mathbb{P}\left\{x+\omega_0 \in E\,\,\text{and} \,\, f_{\omega}\left(I\right)  \subset  \Delta\right\}>q.
\end{align*}
Note that, due to the expansion in $E$, the event $\left\{x+\omega_0 \in E\,\,\text{and} \,\, f_{\omega}\left(I\right)  \subset  \Delta\right\}$ is contained in the event
\begin{equation}\label{ei}
E(I):= \left\{f_{\omega}(I) \subset \Delta,\, |f_{\omega}(I)|\ge \frac{\eta}{4}\,\text{and} \, \left|df_{\omega}\vert_{I} \right|>1 \right\}
\end{equation}
Hence
\begin{equation}\label{chillbachata}
\mathbb{P}\left(E(I) \right)>q.
\end{equation}
By Lemma \ref{hacecalor} with $\eta_1=\frac{\eta}{4}$, we find $K_1$ and $\iota_1$ such that 
\begin{equation}\label{laplaya}
\inf_{\{I \subset \Delta,|I| \ge  \frac{\eta}{4}\}} \mathbb{P}\left\{I \,\,\text{is} \,\, (\omega,K_1,R^{K_1})\text{-full branch expanding}\right\} > \iota_1.
\end{equation}
As a consequence, due to \eqref{ei}
\begin{align*}
&\mathbb{P}\left\{I\,\text{is}\, (\omega,K_1+1,R^{K_1})\,\text{full-branch expanding}  \ \right\} \ge \\ &\int_{E(I)} \mathbb{P}\left
\{f(I+\omega_0)\,\text{is}\,(\theta(\omega),K_1,R^{K_1}) \text{full-branch expanding}\right\}d\mathbb{P}(\omega_0) \ge q\iota_1,
\end{align*}
 which concludes the proof. 
\end{TheoremproofE}
\section{Large deviations estimates for slow recurrence to the critical set and for the Lyapunov exponent}
\subsection{Slow recurrence to the critical set}

In this subsection we show that $(\sigma,R)$-predominantly expanding  RDSs satisfy  Hypothesis \ref{(M)}.
Let $f^n_{\omega}$ be a $(\sigma,R)$-predominantly expanding random dynamical system with $f$ as a generating map. The singular set $\mathcal{SC}$ of $f$ has the form
\begin{align*}
\mathcal{SC}:= \left\{a_0,a_1,\dots,a_{m-1} \right\},
\end{align*}
where $m=|\mathcal{SC}|$ denotes the cardinality of the singular set. 
 Using the same notation as in Hypothesis \ref{(M)}, let 
\begin{align*}
Z_n(\delta,\omega,x)&:= \sum_{i=0}^{n-1}- \log\left(\mathrm{dist}_{\delta}\left(f^i_{\omega}(x),\mathcal{SC}_{\theta^i(\omega)}\right)\right).
\end{align*}
By Birkhoff ergodic Theorem, we see that 
\begin{align*}
\lim_{n \to \infty} \frac{Z_n(\omega,x,\delta)}{n} = \int_{\mathbb{S}^1\times \Omega} \log\left(\mathrm{dist}_{\delta}\left(x,\mathcal{SC}_{\omega}\right)\right)  d\mu(x)d\mathbb{P}(\omega) \qquad \mathbb{P}\otimes \mu \,\, \text{a.s.}
\end{align*}
The main result of this subsection is  the following 
\begin{teorema}\label{largedeviazioni}
Let $f^n_{\omega}$ be a $(\sigma,R)$-predominantly expanding RDS  and let  
\begin{align*}
H(\delta):=\sqrt{\delta}\left(1+\log\left(\frac{1}{\delta}\right)\right).
\end{align*}
Then, there exists $p \in (0,1)$ such that
\begin{align*}
  \mathbb{P}_x( Z_n(\omega,x,\delta)> n H(\delta) ) \le p^n \qquad \forall n \ge 0,
\end{align*}
for $\delta$ small enough.
In particular, Hypothesis \ref{(M)} is fulfilled. 
\end{teorema}
\begin{proof}
By the Markov Inequality, for $1<\zeta< e$, we have
\begin{align*}
  \mathbb{P}\left( Z_n(\omega,x,\delta)> n H(\delta) \right) \le \frac{\mathbb{E}\left[\prod_i \zeta^{-\log\left(\mathrm{dist}_{\delta}\left(f^i_{\omega}(x),\mathcal{SC}_{\theta^i(\omega)}\right)\right)}\right]}{\zeta^{H(\delta)n}}
\end{align*}
In order to simplify the notation, we write $ b_i(\omega,x,\delta):=-\log\left(\mathrm{dist}_{\delta}\left(f^i_{\omega}(x),\mathcal{SC}_{\theta^i(\omega)}\right)\right)$. Note that 
\begin{equation}\label{induction}
\mathbb{E}\left[ \prod_{i\le n} \zeta^{b_i(\omega,x,\delta)}  \right] = \mathbb{E}\left[\prod_{i\le n-1}\zeta^{b_i(\omega,x,\delta)}\mathbb{E}\left[\zeta^{b_n(\omega,x,\delta)}|\mathcal{F}_{n-1}\right](\omega)\right]
\end{equation}
Due to the the identity $\mathcal{SC}_{\omega}=\mathcal{SC}-\omega_0$ and the fact that $\mathcal{SC}_{\theta^n(\omega)}$ depends explicitly on $\omega_n$, whilst $f^n_{\omega}(x)$ depends on $\omega_0,\dots,\omega_{n-1}$, we have,  
\begin{align*}
\mathbb{E}\left[\zeta^{b_n(\omega,x,\delta)}|\mathcal{F}_{n-1}\right](\omega) &= \int_{(-\sigma,\sigma)} \zeta^{-\log\left(\mathrm{dist}_{\delta}(f^n_{\omega}(x)+\omega_n,\mathcal{SC})\right)} d\mathbb{P}(\omega_n)\\
&= \frac{1}{2\sigma} \int_{\mathbb{S}^1} \zeta^{-\log\left(\mathrm{dist}_{\delta}(y,\mathcal{SC})\right)} 1_{B (f^{n}_{\omega}(x),\sigma )}(y)dy \\
&=     \frac{1}{2\sigma}  \int_{\mathbb{S}^1} 1_{B(f^n_{\omega}(x),\sigma)}(y) \left(\frac{1}{\mathrm{dist}_{\delta}(y,\mathcal{SC})}\right)^{\log(\zeta)} \d y  \\
&= \frac{1}{2\sigma}\int_{\mathbb{S}^1 \setminus B(\mathcal{SC},\delta)}1_{B (f^n_{\omega}(x),\sigma)}dy  + \frac{1}{2\sigma}\int_{B(\mathcal{SC},\delta)}1_{B(f^{n}_{\omega}(x),\sigma)}(y)\left(\frac{1}{\mathrm{dist}_{\delta}(y,\mathcal{SC})}\right)^{\log(\zeta)}d(y) 
\\
&\le  1 + \frac{2|\mathcal{SC}|}{\sigma } \int_{(0,\delta)}\left(\frac{1}{x}\right)^{\log(\zeta)} dx \\
&=   1 + \frac{2|\mathcal{SC}|}{\sigma}\frac{\delta^{1-\log(\zeta)}}{1-\log(\zeta)}.
\end{align*}
Substituting this estimate into \eqref{induction}, we achieve by induction in $n$ that
\begin{equation}\label{bondo3}
\mathbb{E}\left[ \prod_{i\le n} \zeta^{b_i(\omega,x,\delta)}  \right] \le \left(\frac{1 + \frac{2|\mathcal{SC}|}{\sigma }\frac{\delta^{1-\log(\zeta)}}{1-\log(\zeta)}}{\zeta^{H(\delta)}}\right)^n.
\end{equation}
Furthermore,
\begin{align*}
\frac{d \zeta^{H(\delta)}}{d\delta} &= \log(\zeta) \zeta^{H(\delta)}\left(\frac{1}{2\sqrt{\delta}} + \frac{\log(\frac{1}{\delta})}{2\sqrt{\delta}}- \frac{1}{\sqrt{\delta}}  \right)\\&= \frac{1}{2\sqrt{\delta}}\log(\zeta) \zeta^{H(\delta)} \left(\log\left(\frac{1}{\delta}\right)-1\right)\\ &> \frac{1}{2\sqrt{\delta}}\log(\zeta) \,\, \text{as} \, \delta \to 0.
\end{align*}
Thus, for sufficiently small $\delta$,
\begin{equation}\label{bondo1}
\zeta^{H(\delta)}-1 > \lim_{t \to 0} \int_t^{\delta} 
 \log(\zeta) \frac{1}{2\sqrt{x}}dx = \log(\zeta) \sqrt{\delta}.
\end{equation}
Choose $\zeta$ sufficiently close to one so that $\frac{1}{2}<1-\log(\zeta)$. Then, if $\delta$ is small enough, we have that
\begin{equation}\label{bondo2}
\log(\zeta)\sqrt{\delta} > \frac{2|\mathcal{SC}|}{\sigma }\frac{\delta^{1-\log(\zeta)}}{1-\log(\zeta)}.
\end{equation}
As a consequence,
\begin{align*}
p := \left(\frac{1 + \frac{2|\mathcal{SC}|}{\sigma}\frac{\delta^{1-\log(\zeta)}}{1-\log(\zeta)}}{\zeta^{H(\delta)}}\right) <1,
\end{align*}
 which, by  \eqref{bondo3},  concludes the proof.
\end{proof}

\subsection{Positivity of Lyapunov exponent and large deviations}
The main result of this subsection is the following
\begin{proposicao}\label{ledev}
Let $f^n_{\omega}$ be $(\sigma,R)$-predominantly expanding RDS. Then, the Lyapunov exponent $\lambda$ is positive and satisfies
\begin{equation}\label{crazy1}
\lambda> Z(h),
\end{equation}
where  $Z(h)$ is given by \eqref{himpl2} in  Definition \ref{admissibleRDS}. 
Furthermore, if $\alpha>0$ is taken small enough so that
\begin{equation}\label{new}
\bar{Z}(h):=  \log(R)(1-h)+\frac{\alpha+1}{2\sigma}\int_{C} \log|df(z)|dz>0,
\end{equation}
where $C$ is the contracting region of $f$, then there exists a constant $D \ge 1$ and $q \in (0,1)$ such that
\begin{equation}\label{crazy2}
\sup_{x \in \mathbb{S}^1} \mathbb{P}\left( \sum_{i=0}^{n-1}\log\left|df_{\theta^i(\omega)}(x)\right|<\bar{Z}(h)n\right) \le Dq^n.
\end{equation}
In particular, Hypothesis \ref{(L)} is fulfilled. 
\end{proposicao}
{\em Proof.}
By Proposition \ref{SoyAquel}, the Markov chain associated to the RDS $f^n_{\omega}$ has an  unique ergodic stationary measure $\mu$ such that $\log\left|df(x+\omega_0)\right| \in L^1(\mathbb{P}\otimes \mu)$. As a consequence
\begin{equation}\label{loteria}
\lambda = \int_{\mathbb{S}^1}\int_{\Omega}\log\left|df(x+\omega_0)\right|d\mu(x)d\mathbb{P}(\omega) \ge \inf_{x \in \mathbb{S}^1}\int_{\Omega}\log\left|df(x+\omega_0)\right|d\mathbb{P}(\omega)
\end{equation}
Furthermore, equations \eqref{himpli1} and \eqref{himpl2} imply that for any $x \in \mathbb{S}^1$ we have
\begin{align*}
\int_{\Omega}\log|df(x+\omega_0)|d\mathbb{P}(\omega)&= \frac{1}{2\sigma}\int_{B(x,\sigma)\cap G}\log|df(z)| dz + \frac{1}{2\sigma}\int_{B(x,\sigma)\setminus G} \log|df(z)|dz \\
&\ge \log(R)(1-h) + \frac{1}{2\sigma}\int_{B(x,\sigma)\setminus G}\log|df(z)| dz>Z(h)>0.
\end{align*}

It remains to establish the large deviations estimate  \eqref{crazy2}.  First, in order to simplify the notation, let $V:=-\int_C\log(|df(x)|)dx>0$.
Then  
\begin{equation}\label{zetauno}
\bar{Z}(h)= (1-h)\log(R)-\frac{\alpha+1}{2\sigma}V.
\end{equation}
For sufficiently large $l \in \mathbb{N}$, we consider the sets
\begin{equation}\label{defPk}
P(k):= \left\{x \in \mathbb{S}^1 \mid R^{-\frac{k}{l}}>|df(x)|>R^{-\frac{(k+1)}{l}} \right\}.
\end{equation}
Then $\{ P(k)\}_{k \ge -l}$ is a partition of $\mathbb{S}^1 \setminus G$ whose diameter is small if $l$ is large  and the sets $\{P(k)\}_{k \ge 0}$ are a partition of $C$  with the same properties. Then, for any sequence $\{x_k\}_{k \ge 0}$ such that $x_k \in P(k)$ for  $l$ large enough we have   
\begin{align*}
 \sum_{k \ge 0} \log|df(x_k)|\mathrm{Leb}(P(k)) \approx -V.
\end{align*}
Furthermore, as a consequence of the definition of $P(k)$ in \eqref{defPk}, we have that  
\begin{align*}
\left|\sum_{k \ge 0} \log|df(x_k)|\mathrm{Leb}(P(k))- \sum_{k \ge 0} \log(R^{-\frac{(k+1)}{l}})\mathrm{Leb}(P(k)) \right| \le \log\left(R^{-\frac{1}{l}}\right)\mathrm{Leb}(C).
\end{align*}
We can choose  $l$ large enough such that
\begin{equation}\label{mirival}
\frac{\sum_{k\ge 0} \frac{(k+1)}{l}\mathrm{Leb}(P(k))}{2\sigma}< \left(\frac{1}{2\sigma\log(R)}+ \frac{\alpha}{4\sigma\log(R)}\right)V.
\end{equation}

Fix $n \in \mathbb{N}$ and $x \in \mathbb{S}^1$. Given $k \ge -l$ let $\tilde{P}_{n,k}(\omega,x)$ be the number of points among $x+\omega_0, \dots, f^{n-1}_{\omega}(x)+\omega_{n-1}$ that belong to $P(k)$:
\begin{equation}\label{100h}
\tilde{P}_{n,k}(\omega,x):= \text{card}\{i \in \{0,\dots, n-1 \} \mid f^i_{\omega}(x) +\omega_i \in P(k)  \},
\end{equation}
 Analogously, we define
\begin{align*}
\tilde{G}_n(\omega,x) &:= \text{card}\{i \in \{0,\dots, n-1 \} \mid f^i_{\omega}(x)+\omega_i \in G  \}.
\end{align*}
Recall that $G$ is the super-expanding region of $f$ and $\bigcup_{k \ge -l}P(k) = \mathbb{S}^1 \setminus G$.
\\
By definition of $G$ and $P(k)$, we have
\begin{equation}\label{lowerb}
S_n(\omega,x):=\sum_{i=0}^{n-1}\log|df_{\theta^i(\omega)}(f^i_{\omega}(x))| \ge \tilde{G}_{n}(\omega,x)\log(R) - \sum_{k =0}^{\infty}\tilde{P}_{n,k}(\omega,x) \frac{(k+1)}{l}\log(R).
\end{equation}
Given $x \in \mathbb{S}^1$, consider the set of noise realizations 
\begin{align*}
A_n(x) := \left\{ \tilde{G}(n)(\omega,x) - \sum_{k \ge 0} \tilde{P}_{n,k}(\omega,x) \frac{(k+1)}{l} < \left(1-h-\left(\frac{1+\alpha}{2\sigma \log(R)}\right)V\right)n \right\}.
\end{align*}
Note that by \eqref{zetauno} and  \eqref{lowerb} we have 
\begin{equation}\label{inclusion}
\left\{ S_n(\omega,x)< \bar{Z}(h)n \right\} \subset A_n(x).
\end{equation}
In order to estimate the measure of $A_n(x)$, we introduce the set $A^1_n(x)$ of all noise realizations for which the orbit of $x$  enters one of the sets  $P(k)$ for  $k \ge n$ at least once:
\begin{align*}
A^1_n(x) &:= \left\{\exists k \ge n \mid \tilde{P}_{n,k}(\omega,x) >0  \right\}.
\end{align*}
 We further consider the set of noise realizations for which the orbit of $x$ visits $G$ with  frequency smaller than  $1-h- \frac{\alpha}{8\sigma\log(R)}V$:
\begin{align*}
A^2_n(x) := \left\{\tilde{G}_{n}(\omega,x)\le \left(1-h-\frac{\alpha }{8\sigma \log(R)}V\right)n \right\}.
\end{align*}
The set $A_n(x)$ can be decomposed as 
\begin{equation}\label{decomposition}
A_n(x) := A^3_n(x) \sqcup \left( A_n(x) \cap \left(A^1_n(x) \cup A^2_n(x)\right)\right),
\end{equation}
where
\begin{equation}\label{name}
\begin{array}{ll}
A^3_n(x) &=  A_n(x) \setminus \left(A^1_n(x)\cup A^2_n(x)\right)\\
&= \left\{ \sum_{k = 0}^{n-1} \tilde{P}_{n,k}(\omega,x) \left(\frac{k+1}{l}\right)  > \tilde{Z}n \right\},
\end{array}
\end{equation}
and
\begin{equation}\label{name2}
\tilde{Z} := \left(\frac{3\alpha}{8\sigma \log(R)}+\frac{1}{2\sigma \log(R)}\right)V.
\end{equation}
As
\begin{align*}
\mathbb{P}(A_n(x)) \le \mathbb{P}\left(A^1_n(x)\right) + \mathbb{P}(A^2_n(x)) + \mathbb{P}(A^3_n(x)),
\end{align*}
 Proposition \ref{ledev}  is to follow from the estimates on the measure of $A^1_n(x)$, $A^2_n(x)$ and $A^3_n(x)$, in  Lemmas \ref{motivation}, \ref{coollemma} and \ref{importandev} below.
\begin{lema}\label{motivation}
There exists a $D_0 \ge 1$, such that, for all $n \ge 0$,
\begin{align*}
\mathbb{P}\left(A^1_n(x)\right) \le D_0hR^{-\frac{n}{2l}}.
\end{align*}
\end{lema}
\begin{proof}
 Let 
\begin{align*}
\tilde{P}(n):= \bigcup_{k \ge n}P(k),
\end{align*}
then
\begin{align*}
\mathbb{P}\left(A^1_n(x)\right) \le \sum_{i=0}^{n-1} \mathbb{P}\left( f^i_{\omega}(x)+x_i \in \tilde{P}(n)  \right) \le \sum_{i=0}^{n-1} \sup_{y \in \mathbb{S}^1}\mathbb{P}\left( y+\omega_i \in \tilde{P}(n)\right).
\end{align*}
By item $(1)$ Definition \ref{admissibleRDS}, the critical set of $f$ is non-degenerated. Then, there exists a constant $C \ge 1$ such that $\mathrm{Leb}(P(k))\le C R^{-\frac{n}{l}}$, for all $n \ge 0$. As a consequence, we have that
\begin{align*}
\mathbb{P}\left( y+\omega_i \in \tilde{P}(n)  \right) \le \frac{\mathrm{Leb}(\tilde{P}(n))}{2\sigma}
&\le D R^{-\frac{n}{l}},
\end{align*}
which implies that 
\begin{align*}
\mathbb{P}\left(A^1_n(x)\right) \le  DnR^{-\frac{n}{l}} \le D_0R^{-\frac{n}{2l}},
\end{align*}
for some $D_0 \ge 1$.
\end{proof}
\begin{lema}\label{coollemma}
 There exists $\tilde{q}_1 \in (0,1)$ such that, for all $n \ge 0$
\begin{align*}
\mathbb{P}\left(A^2_n(x)\right) \le \tilde{q}_1^n.
\end{align*}
\end{lema}
\begin{proof}
Observe that $A^2_n(x)$ can be rewritten as
\begin{align*}
A^2_{n}(x)= \left\{ \sum_{i=0}^{n-1}1_{\mathbb{S}^1 \setminus G}(f^i_{\omega}(x)+\omega_i)>\left(h+ \frac{\alpha}{4}V\right)n \right\}.
\end{align*}
 By the Markov inequality, for any $\zeta>1$, we obtain that 
\begin{align*}
 \mathbb{P}\left(A^2_{n}(x)\right) \le \frac{\mathbb{E}\left[\zeta^{\sum_{i=0}^{n-1}1_{\mathbb{S}^1\setminus G}\left(f^i_{\omega}(x)+\omega_i\right)}\right]}{\zeta^{\left(h+ \frac{\alpha}{4}V\right)n}}
\end{align*}
We may rewrite the upper term as 
\begin{align*}
\mathbb{E}\left[\zeta^{\sum_{i=0}^{n-1}1_{\mathbb{S}^1\setminus G}(f^i_{\omega}(x)+\omega_i)}\right]  &= \mathbb{E}\left[\prod_{i=0}^{n-1} {\zeta^{\sum_{i=0}^{n-1}1_{\mathbb{S}^1\setminus G}(f^i_{\omega}(x)+\omega_i)}}\right]\\ &= \int_{\omega_0,\dots,\omega_{n-1}} \prod_{i=0}^{n-1} \zeta^{1_{\mathbb{S}^1\setminus G}(f^i_{\omega}(x)+\omega_i)} d\mathbb{P}(\omega_0,\dots\omega_{n-1}) \\
&= \int_{\omega_0,\dots,\omega_{n-2}}\prod_{i=0}^{n-2} \zeta^{1_{\mathbb{S}^1\setminus G}(f^i_{\omega}(x)+\omega_i)}\left[ \int_{\omega_{n-1}} \zeta^{1_{\mathbb{S}^1\setminus G}(f^{n-1}_{\omega}(x)+\omega_{n-1})}d\mathbb{P}(\omega_{n-1})\right] d\mathbb{P}(\omega_0,\dots,\omega_{n-2}). 
\end{align*}
Now observe that, given $\omega_0,\dots,\omega_{n-2}$ fixed, we have that
\begin{align*}
\int_{\omega_{n-1}} \zeta^{1_{\mathbb{S}^1\setminus G}(f^{n-1}_{\omega}(x)+\omega_{n-1})}d\mathbb{P}(\omega_{n-1}) = \zeta\frac{\mathrm{Leb}(\omega_n \mid f^{n-1}_{\omega}(x)+\omega_{n-1} \notin G)}{2\sigma}+ 1- \frac{\mathrm{Leb}(\omega_n \mid f^{n-1}_{\omega}(x)+\omega_{n-1} \notin G)}{2\sigma}.
\end{align*}
Furthermore, by \eqref{himpli1} we have 
\begin{align*}
\int_{\omega_{n-1}} \zeta^{1_{\mathbb{S}^1\setminus G}(f^{n-1}_{\omega}(x)+\omega_{n-1})}d\mathbb{P}(\omega_{n-1}) \le \zeta h + 1-h.
\end{align*}
Repeating this argument recursively we obtain
\begin{align*}
\mathbb{E}\left[\zeta^{\sum_{i=0}^{n-1}1_{\mathbb{S}^1\setminus G}(f^i_{\omega}(x)+\omega_i)}\right]  \le (\zeta h + 1-h)^n.
\end{align*}
Hence, for any $\zeta>1$ 
\begin{align*}
 \mathbb{P}\left(A^2_n(x)\right)\le \left(\frac{h\zeta+ 1-h}{\zeta^{h+\frac{\alpha}{4}V}}\right)^n.
\end{align*}
Taking $\tilde{q}_1=\left(\frac{\zeta h+ 1-h}{\zeta^{h+\frac{\alpha}{4}V}}\right)$  for $\zeta$ close  to $1$ concludes the proof. 
\end{proof}
\begin{lema}\label{importandev}
There exists  $\tilde{q}_2 \in (0,1)$ such that, for all $n \ge 0$
\begin{align*}\mathbb{P}\left(A^3_n(x)\right) \le \tilde{q}_2^n.
\end{align*}
\end{lema}
\begin{proof}
By the Markov inequality and \eqref{name}, for all $\zeta>1$.
\begin{equation}\label{Chebishev}
\mathbb{P}(A^3_n(x)) \le \frac{\mathbb{E}\left[\zeta^{ \sum_{k = 0}^{n-1} \tilde{P}_{n,k}(\omega,x)\frac{k+1}{l}}\right]}{\zeta^{n\tilde{Z}}}.
\end{equation}
Define  for $j = 0,\dots, n-1$
\begin{align*}
M_j(\omega,x) := \prod_{k=0}^{n-1}\zeta^{\frac{k+1}{l}1_{P(k)}(f^j_{\omega}(x)+\omega_j)},
\end{align*}
and let
\begin{align*}
\pi_j(A) = \mathbb{P}((\omega_0,\dots,\omega_{j-1}) \in A) \qquad \forall A \subset \mathcal{B}([-\sigma,\sigma]^j).
\end{align*}
Then, by \eqref{100h}
\begin{align*}
\mathbb{E}\left[\zeta^{ \sum_{k = 0}^{n-1} \tilde{P}_{n,k}(\omega,x)(k+1)}\right] &= \mathbb{E}\left[\prod_{k=0}^{n-1}\zeta^{\tilde{P}_{n,k}(\omega,x)(k+1)}   \right]
= \mathbb{E}\left[\prod_{k=0}^{n-1}\prod_{j=0}^{n-1}\zeta^{(k+1)1_{P(k)}(f^j_{\omega}(x)+\omega_j)} \right] \\
&=\mathbb{E}\left[\prod_{j=0}^{n-1}\left[\prod_{k=0}^{n-1} \zeta^{(k+1)1_{P(k)}(f^j_{\omega}(x)+\omega_j)} \right]\right] 
= \int_{\Omega}  \prod_{j=0}^{n-1} M_j(\omega,x)d\pi_n.
\end{align*}
By Fubini's theorem
\begin{equation}\label{fubini}
\int_{\Omega} \prod_{j=0}^{n-1}M_j(\omega,x)d\pi_{n}= \int_{[-\sigma,\sigma]^{n-1}}\prod_{j=0}^{n-2}{M_j(\omega,x)}\left[\int_{[-\sigma,\sigma]}M_{n-1}(\omega,x)d\mathbb{P}(\omega_n)\right]d\pi_{n-1}.
\end{equation}
For given $(\omega_0,\dots,\omega_{n-2}) \in [-\sigma,\sigma]^{n-1}$, we may rewrite
\begin{align*}
\int_{[-\sigma,\sigma]}M_{n-1}(\omega,x)d\mathbb{P}(\omega_{n-1}) &= \int_{[-\sigma,\sigma]} \prod_{k=0}^{n-1} \zeta^{\frac{k+1}{l}1_{P(k)}(f^{n-1}_{\omega}(x)+\omega_{n-1})}d\mathbb{P}(\omega_{n-1})\\
&\le \sup_{y \in \mathbb{S}^1} \int_{[-\sigma,\sigma]} \prod_{k=0}^{n-1} \zeta^{\frac{k+1}{l}1_{P(k)}(y+\omega_{n-1})}\mathbb{P}(\omega_{n-1})
\\
&\le \sup_{y \in \mathbb{S}^1} 1- \frac{\mathrm{Leb}\left(\cup_{k=0}^{n-1 }P(k) \cap B_{\sigma}(y)\right)}{2\sigma}+ \sum_{k=0}^{n-1} \zeta^{\frac{(k+1)}{l}}\frac{\mathrm{Leb}\left(P(k) \cap B_{\sigma}(y)\right)}{2\sigma} \\
&\le \sup_{y \in \mathbb{S}^1} 1+ \sum_{k=0}^{n-1}  (\zeta^{\frac{(k+1)}{l}}-1)\frac{\mathrm{Leb}\left(P(k) \cap B_{\sigma}(y)\right)}{2\sigma} \\
&\le 1+ \sum_{k=0}^{n-1}  (\zeta^{\frac{(k+1)}{l}}-1)\frac{\mathrm{Leb}\left(P(k) \right)}{2\sigma}
\end{align*} 
As a consequence 
\begin{equation}\label{nickchata}
\int_{\Omega} \prod_{j=0}^{n-1}M_j(\omega,x)d\pi_{n-1}\le  W(\zeta) \int_{\Omega} \prod_{j=0}^{n-1}M_j(\omega,x)d\pi_{n-2},
\end{equation}
where 
\begin{align*}
W(\zeta)= 1+ \sum_{k=0}^{n-1} (\zeta^{\frac{(k+1)}{l}}-1)\frac{\mathrm{Leb}\left(P(k)\right)}{2\sigma}.
\end{align*}
Repeating this argument, we obtain
\begin{equation}\label{exponential}
\mathbb{P}(A^3_n(x)) \le \left(\frac{W(\zeta)}{\zeta^{\tilde{Z}}}\right)^n. 
\end{equation}
Note that $W(1)=1$ and, by \eqref{mirival} and \eqref{name2}
\begin{align*}
dW(1)&= \sum_{k=0}^{n-1}
\frac{(k+1)}{l}\frac{\mathrm{Leb}\left(P(k)\right)}{2\sigma} < \left(\frac{1}{2\sigma \log(R)}+\frac{\alpha}{4\sigma \log(R)}\right)V< \tilde{Z}.
\end{align*}
This implies that  $W(\zeta)$ grows slower than $\zeta^{\tilde{Z}}$ near $\zeta=1$. Therefore 
\begin{align*}
q:= \left(\frac{W(\zeta)}{\zeta^{\tilde{Z}}}\right)<1,
\end{align*}
for $\zeta>1$ sufficiently close to $1$,  which concludes the  proof.
\end{proof}

\begin{remark}
To facilitate future applications of the above result, it is important to highlight the proof of Proposition \ref{crazy1} requires only the use of Definition \ref{admissibleRDS}, items $(1),$ $(2)$ and $(3)$.
\end{remark}

\section{Proof of Theorem \ref{Teorema}}
In this final section, we prove Theorem \ref{Teorema}. It follows from Theorem \ref{teb} and the following proposition.
\begin{proposicao}\label{example2}
Let $f^n_{\omega}$ be the random dynamical system generated by  random composition of the maps in \eqref{example}, i.e. maps of the form
\begin{align*}
f_{\omega}(x):= L\sin(2\pi (x+\omega_0)) \, \,(\mathrm{mod}\,1).
\end{align*}
Suppose $L \ge 3$. Then,   $f^n_{\omega}$ is $(\sigma,R)$-predominantly expanding for $R:= \max\{3,L^{\frac{1}{2}} \}$ and  $\sigma > \frac{1}{R}+ c_1\frac{R}{\pi^2L}$,
with $c_1:= \frac{4\pi}{\sqrt{16 \pi^2-1 }}$.
\end{proposicao}
\begin{proof}
Recall that $f^n_{\omega}$ is $(\sigma,R)$-predominantly expanding if and only if the following two conditions are satisfied:
\begin{itemize}
    \item[(1)] the circle map
\begin{align*}
f(x):=  L\sin(2\pi x)\,\,  (\mathrm{mod}\,1)
\end{align*}
satisfies the conditions stated in  Definition \ref{admissibleRDS} $(1)$;
\item[(2)] 
the inequalities in \eqref{reF} and \eqref{noiseRange} are satisfied.
\end{itemize}
Item $(1)$ follows from the definition of $f$.
 To verify $(2)$, we check  the inequalities in \eqref{reF} and \eqref{noiseRange}.
 Equation \eqref{reF} follows from the inequality 
\begin{align*}
 -\frac{c_1}{\pi^2 L} \le \int_{C}\log|df(x)|dx,
\end{align*}
whilst \eqref{noiseRange} follows from \eqref{dierre} and the fact that
\begin{align*}
\mathrm{Leb}\left(\mathbb{S}^1 \setminus G\right)&= \mathrm{Leb}\left\{x \in \mathbb{S}^1 \mid |2\pi L \cos(2\pi x)|<R \right\} \\
&\le \mathrm{Leb}\left\{\bigcup_{\left\{x_0= \frac{\pi}{4},\frac{3\pi}{4}\right\}}\left\{x \in \mathbb{S}^1 \mid \text{dist}(x,x_0)<c_1 \frac{R}{4\pi^2 L} \right\} \right\}\\
&= c_1 \frac{R}{\pi^2 L}.
\end{align*}
This completes the proof.

\end{proof}
\section*{Acknowledgments}
 JL and GT are supported by the EPSRC Centre for Doctoral Training in Mathematics of Random Systems: Analysis, Modelling and Simulation (EP/S023925/1).
\appendix

\bibliographystyle{plain} 
\bibliography{mybib}
\end{document}